\documentclass[twocolumn]{autart}    

\usepackage{graphicx}          
\usepackage{epstopdf}
\usepackage{fancyhdr}
\usepackage{indentfirst}
\usepackage{amsmath}
\usepackage{amssymb}
\usepackage{natbib}
\usepackage{color}
\usepackage{latexsym}
\usepackage{amsfonts}
\usepackage{amssymb,amsmath,amsfonts}
\usepackage{amsmath,mathrsfs,bm,url,times}
\usepackage{caption}
\usepackage{subcaption}
\usepackage{textcomp}
\usepackage{epsfig}
\usepackage{tikz}
\usepackage{tkz-berge}
\usepackage{lipsum}
\usetikzlibrary{calc}
\usetikzlibrary{shapes,fit}
\usepackage{verbatim}
\usepackage{setspace}
\usetikzlibrary{shapes,arrows}
\newtheorem{proposition}{Proposition}
\newtheorem{theorem}{Theorem}
\newtheorem{lemma}{Lemma}

\newtheorem{assumption}{Assumption}

\newtheorem{remark}{Remark}

\usepackage[toc, page]{appendix}
\usepackage{float} 
\usepackage{enumitem}
\usepackage{url}
\usepackage{adjustbox}
\usepackage{booktabs}
\usepackage{makecell}
\usepackage{pifont}
\usepackage{threeparttable} 
\usepackage{makecell}
\usepackage{bbding}   

\usepackage{algorithm}
\usepackage{algorithmic}

\usepackage{array}
\newcolumntype{M}[1]{>{\centering\arraybackslash}m{#1}}
\newcolumntype{N}{@{}m{0pt}@{}}
\usepackage{multirow}

\usepackage{stfloats}
\usepackage{caption}

\usepackage{xcolor, colortbl}
\definecolor{LightGray}{gray}{0.9}
\definecolor{LightCyan}{rgb}{0.88,1,0.8}

\usepackage{color}
\definecolor{gray}{RGB}{128,128,128}

\newcommand{\LX}[1]{{\color{black} #1}}

\newcommand{\blue}[1]{{\color{black} #1}}
\newcommand{\red}[1]{{\color{black} #1}}

\newcommand{\xl}[1]{{\color{black} #1}}
\newcommand{\xul}[1]{{\color{blue} #1}}

\usepackage{hyperref}
\hypersetup{colorlinks,citecolor=blue,linkcolor=blue,urlcolor=blue}
\allowdisplaybreaks
\begin{document}

\begin{frontmatter}

\title{Quantized Distributed Nonconvex Optimization Algorithms\\ with Linear Convergence \LX{under the Polyak--${\L}$ojasiewicz Condition}} 

\thanks[footnoteinfo]{This work was supported by the National Natural Science Foundation of China 62133003, 61991403 \& 61991400, the Knut and Alice Wallenberg Foundation, the
Swedish Foundation for Strategic Research, \LX{and} the Swedish Research Council. \blue{A preliminary version of this paper has been accepted} \LX{for} the 61th IEEE Conference on Decision and Control, Cancun, Mexico, December 6-9, 2022. (Corresponding author:
Tao~Yang.)}

\author[neu,uvic]{Lei Xu}\ead{2010345@stu.neu.edu.cn},
\author[kth]{Xinlei Yi}\ead{xinleiy@kth.se},   
\author[neu]{Jiayue Sun}\ead{\\sunjiayue@ise.neu.edu.cn},             
\author[uvic]{Yang Shi}\ead{\\yshi@uvic.ca},
\author[kth]{Karl H. Johansson}\ead{\\kallej@kth.se},
\author[neu]{Tao Yang}\ead{yangtao@mail.neu.edu.cn}

\address[neu]{State Key Laboratory of Synthetical Automation for Process Industries, Northeastern University, 110819, Shenyang, China}
\address[kth]{School of Electrical Engineering and Computer Science, KTH Royal Institute of Technology,\\ and Digital Futures, 10044, Stockholm, Sweden}  
\address[uvic]{Department of Mechanical Engineering, University of Victoria, Victoria, BC V8W 2Y2, Canada}             


\begin{abstract}                          
This paper considers distributed optimization for minimizing the average of local \blue{nonconvex} cost functions, by using local information exchange over undirected communication networks.
\LX{To reduce the required communication capacity, we introduce an encoder--decoder scheme}.
By integrating them with distributed gradient tracking and proportional integral algorithms\LX{,} \blue{respectively}, we then propose two quantized distributed nonconvex optimization algorithms.
Assuming the global cost function satisfies the Polyak--{\L}ojasiewicz condition, which does not require the global cost function to be convex and the global minimizer is not necessarily unique, we show that our proposed algorithms linearly converge to a global optimal point and \LX{that larger quantization level leads to faster convergence speed}.
Moreover, \blue{we show that} a low data rate is \blue{sufficient to guarantee} linear convergence when the \LX{algorithm parameters} are properly chosen.
The theoretical results are illustrated by \LX{numerical examples}.
\end{abstract}


\end{frontmatter}

\section{Introduction}
Distributed optimization, which can be traced back to \cite{Tsitsiklis1986TAC,Bertsekas1989parallel}, has received a growing and renewed interest over the last decade due to its wide applications in resource allocation, machine learning, and sensor networks, just to name a few.
Various distributed optimization algorithms have been developed. \LX{For an overview, see} recent survey papers \cite{nedic2018distributed,Yang2019ARC}.
The basic convergence results of distributed optimization algorithms guarantee {\em \LX{sublinear}} convergence to the optimal point \LX{when} the local cost functions are convex, e.g., \cite{Johansson2008CDC,Nedic2009TAC,Zhu2011TAC,Nedic2015TAC}.
When the local cost functions are strongly convex and smooth, linear convergence results are established \LX{(\cite{lu2012zero,Kia2015Aut,wang2010AAC,Xu2015CDC,Shi2015SIAM})}.

Distributed optimization algorithms require the agents to communicate with each other.
Since communication channels have limited bandwidth, distributed optimization algorithms with quantized communications have been developed.
\LX{Early works on quantization for single-agent systems are given in \cite{Larson1967TAC,Curry1969TAC}.
\blue{Extensions have been proposed to deal with the distributed consensus problem with limited communication data rate, e.g.}, \cite{Cai2012Aut,Carli2010Aut,Rikos2020TAC,Lee2020TAC}.
\red{Recent research has focused on \blue{quantized} distributed optimization.}}
For the convex case, \cite{Rabbat2005JSAC,Nedic2008CDC} proposed a quantized distributed incremental and subgradient algorithm, respectively.
These algorithms {\em \LX{sublinearly}} converge to a neighborhood around the optimal point.
\cite{Pu2016TAC} developed a quantized distributed accelerated gradient algorithm and established linear convergence to a neighborhood around the optimal point.
\blue{\cite{Yuan2012distributed} proposed a distributed dual averaging method with quantized communication \red{by using} \LX{a} probabilistic quantizer, \red{and} demonstrated that the proposed \LX{algorithm} {\em \LX{sublinearly}} converges to a global optimal point in expectation.}

Recently, focusing on the strongly convex case, \LX{a} few studies proposed quantized distributed algorithms \LX{that} converge to the exact optimal point.
For example, \cite{Yi2014TCNS} designed a quantized algorithm by integrating the distributed subgradient algorithm and the uniform quantization, while \cite{Doan2020TAC,Zhang2018TAC} developed a quantized gradient algorithm by using the random quantizer and the sign of \LX{the} relative state.
\cite{Xiong2022IS} proposed a quantized mirror descent algorithm \LX{using} time-varying quantizers.
\LX{These algorithms, however, only have {\em \LX{sublinear}} convergence rates.}
\cite{Ma2021Springer,Xiong2021Arxiv} proposed quantized algorithms by equipping the distributed gradient tracking algorithm with uniform quantizers, and established linear convergence to the exact global optimal point for undirected and directed graphs, respectively.
\begin{table*}\label{table1}
\centering
\caption{Comparison of Different Quantized Distributed Algorithms}
\vspace{3mm}
\resizebox{\textwidth}{!}{
\begin{tabular}{lllll} 
\toprule

  \quad\quad\quad\quad\quad\quad\quad \LX{Existing Results in References} & Exact solution & Linear convergence & \makecell[c]{Nonconvex \\ (Local cost functions)}\\
  \midrule
  \makecell[c]{\cite{Rabbat2005JSAC,Nedic2008CDC,Liu2016TCNS}} & \quad\quad~ No & \quad\quad\quad No  & \quad\quad\quad\quad No\\
  \midrule
  \makecell[c]{\cite{Yuan2012distributed,Yi2014TCNS,Zhang2018TAC}\\ \cite{Reisizadeh2019TSP,Doan2020TAC,Xiong2022IS}} & ~~~~~~~ Yes & \quad\quad\quad No  & \quad\quad\quad\quad No\\
  \midrule
  \makecell[c]{\cite{Pu2016TAC,Kajiyama2021TAC}} & \quad\quad~ No & \quad\quad~~ Yes  & \quad\quad\quad\quad No\\
  \midrule
  \makecell[c]{\cite{Ma2021Springer, Xiong2021Arxiv,Lei2020SIAM}} & ~~~~~~~ Yes& \quad\quad~~~Yes & \quad\quad\quad\quad No\\
  \midrule
  ~~\quad\quad\quad\quad\quad\quad\quad\quad\quad This paper & ~~~~~~~ Yes & \quad\quad~~~Yes  & \quad\quad\quad~~ Yes\\
  \bottomrule
  \end{tabular}
}
  \label{tbl:table1}
\end{table*}

\LX{Note that the aforementioned distributed algorithms with linear convergence to the exact optimal point only focus on strongly convex local cost functions.
\xul{However, as demonstrated in \citet[Lemma~2]{Fazel2018PMLR}, the cost function of the LQR problem in reinforcement learning is quadratic and satisfies the Polyak--{\L}ojasiewicz (P--{\L}) condition under appropriate system dynamics and control laws. Moreover, in deep learning, \citet[Theorem~4]{Liu2022ACHA} demonstrated that certain wide neural networks satisfy the P--{\L} condition.}
\xl{This motivates us to consider the P--{\L} condition case.}}

\xl{The main contributions are summarized as follows.

(i) We propose quantized distributed algorithms by integrating the encoder--decoder scheme and the uniform quantizer with the distributed gradient tracking algorithm and distributed proportional integral algorithm, respectively.

(ii)
Assuming that the global cost function satisfies the P--{\L} condition, Theorems~\ref{thm-3} and \ref{thm-1} show that the proposed algorithms linearly converge to an exact global optimal point when the quantization level is larger than a certain threshold.
This is more general than the existing results in \cite{Kajiyama2021TAC,Ma2021Springer,Xiong2022TAC,Lei2020SIAM}, which require that each local cost function is strongly convex.
Table~\ref{tbl:table1} summarizes the comparison between this paper and related studies.

(iii) Moreover, Theorems~\ref{thm-4} and \ref{thm-2} show that the proposed algorithms with arbitrary quantization level can still converge linearly to an exact global optimal point provided that the algorithm parameters are properly chosen.
It is worth noting that \cite{zhao2022NIPS} proposed the BEER algorithm for distributed nonconvex optimization under compressed communication. Unlike the randomized compression operator used in \cite{zhao2022NIPS}, our proposed algorithms employ a uniform quantizer. Theorems~2 and 4 demonstrate linear convergence even for 1-bit data rate, which is the most communication-efficient.}


The remainder of the paper is organized as follows. Section~\ref{sec-2} presents the problem formulation.
\LX{Section~\ref{sec-3.1} introduces an encoder--decoder scheme for quantized communication.}
In Section~\ref{sec-3} and Section~\ref{sec-4}, we propose the distributed gradient tracking and distributed proportional integral algorithms with finite data rates and analyze their performance, respectively.
Section~\ref{sec-5} presents numerical simulation examples.
Finally, concluding remarks are offered in Section~\ref{sec-6}.
To improve the readability, all the proofs are given in the appendix \footnote{\LX{Some of the results related to quantized distributed proportional integral algorithm appeared in the conference paper \cite{Xu2022CDC}. Further to \cite{Xu2022CDC}, this paper significantly adds new contents and analysis. This paper proposes not only the quantized distributed proportional integral algorithm but also the quantized distributed gradient tracking algorithm. Moreover, it contains detailed proofs omitted from the conference version.}}.

\textbf{Notation:} Let $\mathbf{1}_n$ (or $\mathbf{0}_n$) be the $n\times 1$ vector with all ones (or zeros), and $\mathbf{I}_n$ be the $n$-dimensional identity matrix.
$\|\cdot\|$ is the Euclidean vector norm or spectral matrix norm. For a column vector $x=(x_{1},\dots,x_{m})$,  $\|x\|_{\infty}=\max_{1\leq i\leq m}|x_{i}|$.
For a positive semi-definite matrix $\mathcal{M}$,
$\rho(\mathcal{M})$ and $\underline{\rho}(\mathcal{M})$ are the spectral radius and the minimum positive eigenvalue of matrix $\mathcal{M}$, respectively.
The minimum integer greater than or equal to $c$ is denoted by $\lceil c \rceil$.
Let ${\rm{diag}}[a_{1},\dots,a_{n}]$ denote a diagonal matrix with the $i$-th diagonal element being $a_{i}$. Given any differentiable function $f$, $\nabla f$ is the gradient of $f$.
$A\otimes B$ represents the Kronecker product of matrices $A$ and $B$. \LX{$A\preceq B$ if all entries of matrix $A-B$ are not greater than zero, and $A\succ0$ if all entries of matrix $A$ that are greater than zero.}

\section{Problem Formulation}\label{sec-2}

Consider a group of $n$ agents distributed over an undirected graph $\mathcal{G}=(\mathcal{V},\mathcal{E},\mathcal{A})$, where $\mathcal{V}=\{1,2,\ldots,n\}$ \LX{is the vertex set} and $\mathcal{E}\subseteq\mathcal{V}\times\mathcal{V}$ is the set of edges.
$(i,j)\in\mathcal{E}$ indicates that agents $i$ and $j$ can communicate with each other,
and $\mathcal{A}=[a_{ij}]\in\mathbb{R}^{n\times n}$ is the adjacency matrix, where $a_{ij}>0$ if $(j,i)\in\mathcal{E}$, otherwise $a_{ij}=0$.
Let $\mathcal{N}_i=\{j\in\mathcal{V}: a_{ij}>0\}$ and $d_i=\sum_{j=1}^{n}a_{ij}$ denote the \LX{neighbor} set and weighted degree of agent $i$, respectively.
The degree matrix is defined as $\mathcal{D}=\operatorname{diag}[d_1,\ldots,d_n]$.
The graph Laplacian matrix is $L=[L_{ij}]=\mathcal{D}-\mathcal{A}$.
A path from agent $i_{1}$ to agent $i_{k}$ is a sequence of agents $\{i_{1},\dots,i_{k}\}$ such that $(i_{j},i_{j+1})\in\mathcal{E}$ for $j=1,\dots,k-1$. An undirected graph is connected if there exists a path between any pair of distinct agents.

Assume that each agent has a local cost function $f_i:\mathbb{R}^m\rightarrow\mathbb{R}$. The objective is to find an optimizer $x^\star$ to minimize the following optimization problem
\begin{equation}\label{MinFunction}
\min_{x \in \mathbb{R}^m}f(x) = \frac{1}{n} \sum_{i=1}^{n}f_{i}(x).
\end{equation}
Throughout this paper, we make the following assumptions.
\begin{assumption}\label{assum-1}
The undirected graph $\mathcal{G}$ is connected.
\end{assumption}
\begin{assumption}\label{assum-2}
Each local cost function $f_i(x)$ is smooth with constant $L_{f}>0$, i.e.,
\begin{equation}\label{assum2}
\|\nabla f_i(x)-\nabla f_i(y)\|\leq L_{f}\|x-y\|,~\forall x,y\in\mathbb{R}^{m}.
\end{equation}
\end{assumption}
\begin{assumption}\label{assum-3.3}
The optimal set $\mathbb{X}^{\star}=\emph{\emph{argmin}}_{x\in\mathbb{R}^{m}}f(x)$ is nonempty and $f^{\star}=\emph{\emph{min}}_{x\in\mathbb{R}^{m}}f(x)>-\infty$.
\end{assumption}
\begin{assumption}\label{assum-3}
The global cost function $f(x)$ satisfies the Polyak--{\L}ojasiewicz condition with constant $\nu>0$, i.e.,
  \begin{equation}\label{Assum3}
  \frac{1}{2}\|\nabla f(x)\|^{2}\geq\nu(f(x)-f^{\star}),~\forall x\in\mathbb{R}^{m}.
  \end{equation}
\end{assumption}

\begin{remark}\label{rem-1}
Assumptions \ref{assum-1}--\ref{assum-3.3} are common in the literature, e.g., \cite{nedic2018distributed,Yang2019ARC}.
\xl{Assumption~\ref{assum-3} is weaker than strong convexity, and the
global minimizer is not necessarily unique, but every stationary
point is a global minimizer.
Every function that is strongly convex satisfies the P--{\L} condition.
In practical applications, such as reinforcement learning and deep learning, the P--{\L} condition  is satisfied \cite{Fazel2018PMLR,Liu2022ACHA}.}
\end{remark}

The objective of this paper is to propose quantized distributed optimization algorithms with linear convergence under the P--{\L} condition.

\LX{\section{Encoder--Decoder Scheme for Quantized Communication}\label{sec-3.1}
In this section, we introduce an encoder--decoder scheme.}

To begin with, let us consider a uniform quantizer $q[a]$ with $2\mathcal{K}+1$ quantization levels (\cite{Gray1998TIT}), i.e., 
\begin{align}\label{QuantizationRule}
q[a]&=
  \left\{\begin{aligned}
  &j,\quad\quad\quad\frac{2j-1}{2}<a \leq\frac{2j+1}{2},~j=0,\cdots,\mathcal{K},\\
  &\mathcal{K},\quad\quad~~\frac{2\mathcal{K}+1}{2}>a,\\
  &-q[a],~~a\leq-\frac{1}{2}.
  \end{aligned}\right.
  \end{align}
\blue{For this $2\mathcal{K}+1$-level quantizer, the communication channel is required to be capable
of transmitting $\left \lceil \log_2(2\mathcal{K}) \right \rceil$ bits.
Next for a vector $l=[l_1, l_2, \ldots, l_m] \in \mathbb{R}^m$, we define $Q[l]=(q[l_{1}],\dots,q[l_{m}])$.
The quantizer $Q[l]$ is not saturated if $\|l\|_{\infty}\leq\mathcal{K}+\frac{1}{2}$. In this case, the quantization error is bounded, i.e.,
  \begin{equation}\label{QuanError}
  \|l-Q[l]\|_{\infty}\leq\frac{1}{2}.
  \end{equation}}

Next, we introduce an encoder--decoder pair \cite{Lei2020SIAM,Ma2021Springer} for agents to communicate with each other.
First, the following encoder scheme is used to quantize the variable to be transmitted.


\floatname{algorithm}{Encoder}
\renewcommand{\thealgorithm}{}
\begin{algorithm}[!h]
  \label{Encoder2}
  \begin{algorithmic}
  \caption{}
  \STATE \xl{For the vector $\mathcal{C}_j(k)\in\mathbb{R}^{m}$ that requires quantization}, agent $j\in\mathcal{V}$ recursively generates \LX{the} $m$-dimensional quantized output ${z_{j}^{\mathcal{C}}(k)}$, and internal state ${b_{j}^{\mathcal{C}}(k)}$ as follows\LX{: for any $k\geq1$,}\\
  \begin{subequations}\label{encoder2}
  \begin{align}
  z_{j}^{\mathcal{C}}(k)&= Q\left[\frac{1}{s(k-1)}(\mathcal{C}_{j}(k)-b_{j}^{\mathcal{C}}(k-1))\right],\label{Encoder2Aa}\\
  b_{j}^{\mathcal{C}}(k)&= s(k-1)z_{j}^{\mathcal{C}}(k)+b_{j}^{\mathcal{C}}(k-1),\label{Encoder2Ba}
  \end{align}
  \end{subequations}
where the initial value $b_{j}^{\mathcal{C}}(0)=\mathbf{0}_m$, $s(k)=s(0)\mu^{k}>0$ is a decreasing sequence used to adaptively adjust the encoder, and $\mu\in(0,1)$ is a positive constant.
\end{algorithmic}
\end{algorithm}


The following decoder scheme is used to recover the variables sent from agent $j$.

\floatname{algorithm}{Decoder}
\renewcommand{\thealgorithm}{}
\begin{algorithm}[!h]
  \label{Decoder2}
  \begin{algorithmic}
  \caption{}
  \STATE When agent $i\in\mathcal{N}_{j}$ receives the quantized data $z_{j}^{\mathcal{C}}(k)$ from agent $j$, \LX{it} recursively generates an estimate $\hat{\mathcal{C}}_{j}(k)$ \LX{of ${\mathcal{C}_{j}(k)}$ by the following rule: for any $k\geq1$,}\\
  \begin{subequations}\label{decoder2}
  \begin{align}
  \hat{\mathcal{C}}_{j}(k)= s(k-1)z_{j}^{\mathcal{C}}(k)+\hat{\mathcal{C}}_{j}(k-1),\label{Decoder2A}
  \end{align}
  \end{subequations}
where the initial value $\hat{\mathcal{C}}_{j}(0)=\mathbf{0}_m$.
\end{algorithmic}
\end{algorithm}
\begin{remark}\label{rem-3}
From the encoder--decoder scheme, we note that $b_{j}^{\mathcal{C}}(k)$ is a predictor, $s(k)$ is used to adjust the prediction error $\mathcal{C}_{j}(k)-b_{j}^{\mathcal{C}}(k-1)$. Moreover, the initial value $s(0)$ requires to be large enough to guarantee that the quantizer is not saturated, which implies the quantization error is bounded. The positive constant $\mu\in(0,1)$ ensures that the agent gradually improves the accuracy of the estimate for the transmitted variables from neighbors.
\end{remark}

\section{Distributed Gradient Tracking Algorithm with Finite Data Rates}\label{sec-3}
In this section, we propose a distributed gradient tracking algorithm with quantized communication by integrating the distributed gradient tracking algorithm \cite{Nedic2017SIAM,Qu2017TCNS} with the encoder--decoder scheme.
More specifically, the quantized distributed gradient tracking algorithm is given in Algorithm~\ref{alg:B}.

\floatname{algorithm}{Algorithm}
\renewcommand{\thealgorithm}{1}
\begin{algorithm}[!h]
  \caption{Quantized Distributed Gradient Tracking Algorithm}
  \label{alg:B}
  \begin{algorithmic}
  \blue{
  \STATE  $\bm{\mathrm{For~each~agent}}$ $i\in\mathcal{V}.$
	\STATE $\bm{\mathrm{Initialization:}}$\\
     $x_{i}(0)\in\mathbb{R}^{m},~u_{i}(0)=\nabla f_{i}(x_{i}(0)),~\hat{x}_{j}(0)=\hat{u}_{j}(0)=\mathbf{0}_m$.
     }
  \STATE  $\bm{\mathrm{for}}$ $k\geq0:$
	\STATE $\bm{\mathrm{Update:}}$
\begin{subequations}\label{Alg2}
\begin{align}
  x_{i}(k+1)&=x_{i}(k)-\beta\sum_{j=1}^{n}L_{ij}\hat{x}_{j}(k)-\delta u_{i}(k),\label{Alg2a}\\
  u_{i}(k+1)&=u_{i}(k)-\beta\sum_{j=1}^{n}L_{ij}\hat{u}_{j}(k)+\nabla f_{i}(x_{i}(k+1))\notag\\
  &\quad-\nabla f_{i}(x_{i}(k)),\label{Alg2b}
  \end{align}
\end{subequations}
where $\beta$ and $\delta$ are gain parameters.\\

\STATE $\bm{\mathrm{Send:}}$\\
The quantized outputs $z_{i}^{x}(k+1)$ and $z_{i}^{u}(k+1)$ generated by encoder~(\ref{encoder2}) to its neighbors.
%
\STATE $\bm{\mathrm{Receive:}}$\\
The quantized outputs $z_{j}^{x}(k+1)$ and $z_{j}^{u}(k+1)$ generated by encoder~(\ref{encoder2}) from its neighbors.
\STATE $\bm{\mathrm{Compute:}}$\\
The variables $\hat{x}_{j}(k+1)$ and $\hat{u}_{j}(k+1)$ generated by decoder~(\ref{Decoder2A}).
%
\end{algorithmic}
\end{algorithm}

Before stating the main convergence results, we provide the following preliminary results.

The following lemma provides a sufficient condition to ensure that a certain linear matrix inequality holds, which plays a crucial role in ensuring convergence of the consensus error, gradient tracking error, optimization error, and the nonsaturation of the uniform quantizer.

\begin{lemma}\label{lem-1}
Suppose that the parameters $\beta$ and $\delta$ satisfy
  \begin{equation}\label{Prop2beta}
  \beta\in(0,\frac{\sqrt{2}}{2\rho(L)}),
  \end{equation}
  \begin{align}\label{Prop2delta}
  &\delta\in(0,\min\{\frac{\sqrt{c_{1}\Theta_{1}}}{2},\frac{1}{4L_{f}},~\frac{2}{\nu},\frac{1}{8+2L_{f}},\frac{\sqrt{c_{1}}}{8L_{f}}\}),
  \end{align}
where
  \begin{align*}
  &c_{1}=\frac{(1-\varrho^{2}-\LX c_{2})(1-\varrho^{2})}{1+\varrho^{2}},~\varrho=\rho(\bm{\mathrm{I}}_{nm}-\beta\bm{L}-\bm{H}),\\
  &\bm{L}=L\otimes\bm{\mathrm{I}}_{m},~\bm{H}=\frac{1}{n}(\bm{1}_{n}\bm{1}^{T}_{n}\otimes\bm{\mathrm{I}}_{m}),\\
  &\LX{c_{2}\in(0,1-\varrho^{2})},~\Theta_{1}=\min\{\frac{c_{1}}{24L_{f}^{2}},\frac{\nu\Theta_{2}}{2L_{f}^{2}}\},~\LX{\Theta_{2}=\frac{c_{1}}{32L_{f}^{2}}}.
  \end{align*}

Then, the following linear matrix inequality holds:
\begin{equation}\label{Phi}
\Phi\Theta\leq(1-\iota)\Theta,
\end{equation}
where $\iota=\min\{\frac{\LX c_{2}}{2},\frac{\delta}{4}\nu\}$, $\Theta=[\Theta_{1},1,\Theta_{2}]^{T}$
and the nonnegative matrix $\Phi$ is given by
\begin{align}\label{PPhi}
\Phi=\left[
       \begin{array}{cccc}
         \chi_{1} & \chi_{2} & 0\\
         \chi_{3} & \chi_{4} & \chi_{5}\\
         \chi_{6} & 0 & \chi_{7}
       \end{array}
     \right],
\end{align}
where
\begin{align*}
&\chi_{1}=(1+\sigma_{1})\varrho^{2},~\chi_{2}=2\delta^{2}(1+\frac{1}{\sigma_{1}}),\\
&\chi_{3}=(1+\frac{1}{\sigma_{1}})8L_{f}^{2}(\beta^{2}\rho^{2}(L)+\frac{2\delta^{2}L_{f}^{2}}{1-2\delta L_{f}}),\\
&\chi_{4}=(1+\sigma_{1})\varrho^{2}+(1+\frac{1}{\sigma_{1}})8L_{f}^{2}\delta^{2},\\
&\chi_{5}=(1+\frac{1}{\sigma_{1}})\frac{\LX{16L_{f}^{2}\delta(2-\delta\nu)}}{1-2\delta L_{f}},\\
&\chi_{6}=\frac{\delta}{2}L_{f}^{2},~\chi_{7}=1-\frac{\delta}{2}\nu,~\sigma_{1}=\frac{1-\varrho^{2}}{2\varrho^{2}}.
\end{align*}
\end{lemma}

\textbf{Proof} The proof is given in Appendix \ref{Appendix-Lem1}.

The following lemma establishes an upper bound for $\|\Phi^{k}\|$.

\begin{lemma}\label{lem-3}
Suppose that the parameters $\beta$ and $\delta$ are given in Lemma~\ref{lem-1}.
Then, the following inequality holds:
\begin{equation}\label{hatrho}
\|\Phi^{k}\|\leq h\rho^{k}(\Phi),
\end{equation}
where $h=\LX{\sqrt{3}}\frac{\max_{1\leq i\leq 3}\zeta_{i}}{\min_{1\leq i\leq 3}\zeta_{i}}$,
$\zeta=[\zeta_{1},\zeta_{2},\zeta_{3}]^{T}$ is an eigenvector of $\Phi$ corresponding to the spectral radius $\rho(\Phi)$.
\end{lemma}

\textbf{Proof} The proof is given in Appendix \ref{Appendix-Lem3}.

Denote
$F(\bm{x})=\sum_{i=1}^{n}f_{i}(x_{i})$,
$\bm{x}(k)=[x^T_{1}(k),\dots,x^T_{n}(k)]^{T}$,
$\bm{u}(k)=[u^T_{1}(k),\dots,u^T_{n}(k)]^{T}$,
$\bar{x}(k)=\frac{1}{n}(\bm{1}_{n}^{T}\otimes\bm{\mathrm{I}}_{m})\bm{x}(k)$,
$\bar{u}(k)=\frac{1}{n}(\bm{1}_{n}^{T}\otimes\bm{\mathrm{I}}_{m})\bm{u}(k),$
$\bar{\bm{x}}(k)=\bm{1}_{n}\otimes\bar{x}(k)$,
and
$\bar{\bm{u}}(k)=\bm{1}_{n}\otimes\bar{u}(k)$.
The following proposition shows that the quantizer in \eqref{Encoder2Aa} is never saturated by appropriately choosing the proposed algorithm's parameters if the quantization level is larger than a certain threshold.

\begin{proposition}\label{prop2}(nonsaturation).
Suppose that \blue{Assumptions~\ref{assum-1}--\ref{assum-3}} hold. \LX{Consider} Algorithm~\ref{alg:B} and parameters \blue{$\delta$ and $\beta$} \LX{as} given in Lemma~\ref{lem-1}.
Then, for any
\begin{equation}\label{K-level2}
\mathcal{K}\geq\max\{\vartheta_{1},\vartheta_{2}\},
\end{equation}
where\blue{
\begin{align*}
\vartheta_{1}&=\sigma_{3}\sqrt{\frac{h\sigma_{2}}{4\mu^{2}(\mu^{2}-\bar{\rho})}}+\frac{(1+2\beta d)}{2\mu}-\frac{1}{2},\\
\vartheta_{2}&=\sigma_{4}\sqrt{\frac{h\sigma_{2}}{4\mu^{2}(\mu^{2}-\bar{\rho})}}+\frac{\sqrt{nm}L_{f}\beta\rho(L)}{2\mu}+\frac{(1+2\beta d)}{2\mu}-\frac{1}{2},
\end{align*}}
and $d$ is the maximum degree in $\mathcal{G}$,
\begin{align*}
\mu&\in(\sqrt{\bar{\rho}},1),~\bar{\rho}=1-\iota,\\
\sigma_{2}&=\blue{nm}\sigma_{5}\sqrt{1+(1+4L_{f}^{2})^{2}},~\sigma_{3}=\sigma_{6}\sqrt{2+\frac{\LX{4}\sigma_{7}}{\delta(1-2 \delta L_{f})}},\\
\sigma_{4}&=\blue{\beta\rho(L)(1+L_{f})+L_{f}\delta(1+\sqrt{\frac{\LX{4}\sigma_{7}}{\delta(1-2 \delta L_{f})}})},\\
\sigma_{5}&=(1+\frac{1}{\sigma_{1}})2\beta^{2}\rho^{2}(L),~\sigma_{6}=\sqrt{3}\max\{\beta\rho(L),\delta\},\\
\sigma_{7}&=\max\{1-\frac{\delta}{2}\nu,\frac{\delta}{2}L_{f}^{2}\},
\end{align*}
the quantizer \blue{in} (\ref{Encoder2Aa}) is never saturated provided that
\begin{align}\label{ScaF2}
&s(0)\geq\max\Bigg\{\frac{2(C_{x}+\delta C_{u})}{2\mathcal{K}+1},~\frac{2\|\nabla F(\bm{x}(0)-\delta\bm{u}(0))\|_{\infty}}{2\mathcal{K}+1},\notag\\
&\quad\quad\quad\quad\quad\quad\sqrt{\frac{4\|\Lambda(0)\|\mu^{2}(\mu^{2}-\bar{\rho})}{\sigma_{2}}}\Bigg\},
\end{align}
where $C_{x}\geq\|\bm{x}(0)\|_{\infty}$, $C_{u}\geq\|\bm{u}(0)\|_{\infty}$, and $\Lambda(0)=[\|\bm{x}(0)-\bar{\bm{x}}(0)\|^{2},~\|\bm{u}(0)-\bar{\bm{u}}(0)\|^{2},~n(f(\bar{x}(0))-f^{\star})]^{T}$.
\end{proposition}
\textbf{Proof:} The proof is given in Appendix \ref{Appendix-D}.

Proposition \ref{prop2} provides a sufficient condition to ensure that the quantizer \LX{ensures} nonsaturation.
We are now ready to present the convergence result.

\begin{theorem}\label{thm-3}
(high data rate). Suppose that Assumptions~\ref{assum-1}--\ref{assum-3} hold.
Let each agent $i\in\mathcal{V}$ run the Algorithm~\ref{alg:B} with the same $\beta$, $\delta$, $\mathcal{K}$, \blue{$\mu$} and $s(0)$ given in Proposition \ref{prop2}.
Then,
\begin{align}\label{LambdakThm-1}
\blue{\|\bm{x}(k)-\bar{\bm{x}}(k)\|^{2}+n(f(\bar{x}(k))-f^{\star})\leq\sigma_{9}\mu^{2k}.}
\end{align}
where \blue{$\sigma_{9}=\frac{h\sigma_{2}s^{2}(0)}{2\mu^{2}(\mu^{2}-\bar{\rho})}$}.
\end{theorem}
\textbf{Proof:} The proof is given in Appendix \ref{Appendix-E}.

Theorem \ref{thm-3} establishes linear convergence of the proposed algorithm provided that the quantization level is larger than a certain threshold given in (\ref{K-level2}).
The following theorem establishes linear convergence result for arbitrarily low data rate, \xl{even bit rate one, and thus is called low data rate theorem.}

\begin{theorem}\label{thm-4}
(low data rate). Suppose that Assumptions~\ref{assum-1}--\ref{assum-3} hold. Let each agent $i\in\mathcal{V}$ run the Algorithm \ref{alg:B} with $(\mu,\beta,\delta)\in\Pi$, where
\begin{align*}
&\Pi=\{\blue{(\mu,\beta,\delta):\mu\in(\sqrt{\bar{\rho}},1)},~\beta\in(0,\frac{\sqrt{2}}{2\rho(L)}),\\
&\quad\quad~~\delta\in(0,\min\{~\frac{\sqrt{c_{1}\Theta_{1}}}{2},~\frac{1}{4L_{f}},~\frac{2}{\nu},~\frac{1}{8+2L_{f}},~\frac{\sqrt{c_{1}}}{8L_{f}}\}),\notag\\
&\quad\quad~~\LX{\vartheta_{1}\leq \mathcal{K},~\vartheta_{2}\leq \mathcal{K}}\}.
\end{align*}
Then, for any $\mathcal{K}\geq1$ and $s(0)$ satisfying (\ref{ScaF2}) in Proposition~\ref{prop2}, \xl{$\Pi$ is nonempty, and}
\begin{align}\label{LambdakThm-2}
\blue{\|\bm{x}(k)-\bar{\bm{x}}(k)\|^{2}+n(f(\bar{x}(k))-f^{\star})\leq\sigma_{9}\mu^{2k}.}
\end{align}
\end{theorem}
\textbf{Proof:} The proof is given in Appendix \ref{Appendix-F}.
\vspace{-0.3cm}
\xul{
\begin{remark}\label{rem-4}
From \eqref{K-level2}, it can be observed that a smaller $\mu$, i.e., a faster convergence speed, leads to a larger $\mathcal{K}$ and thus a larger quantization level, i.e., a larger communication bandwidth requirement.
Moreover, a smaller $\mu$ requires fewer iterations to achieve a certain level of optimization accuracy.
Note that in Theorem~2, linear convergence of Algorithm~1 is established even for 1-bit data rate, i.e., $\mathcal{K}=1$, which is the most communication efficient.
\end{remark}
}

\vspace{-0.3cm}
\xl{
\begin{remark}\label{rem-5}
For the strongly convex case, \cite{Kajiyama2021TAC,Ma2021Springer,Xiong2022TAC} proposed quantized distributed algorithms with linear convergence.
However, their analysis cannot be used for the P--{\L} condition.
For example, the linear system of inequalities in the above studies use $\|\bar{x}(k)-x^{\star}\|$, where $x^{\star}$ is the unique optimal solution which exists due to the strong convexity.
In our case, the optimal solution is not unique due to the P--{\L} condition.
Therefore, we use $n(f(\bar{x}(k))-f^{\star})$, where $f^{\star}$ is the unique optimal value.
Moreover, as shown in \cite[Lemma~4.1]{Ma2021Springer} and \cite[Lemma~8]{Xiong2022TAC}, the authors leveraged the strong convexity condition to directly apply \cite[Lemma~10]{Qu2017TCNS} for analyzing the upper bound of $\|\bar{x}(k)-x^{\star}\|$.
However, this cannot be used for our case due to the lack of the strong convexity.
Instead, we use the P--{\L} condition to analyze the upper bound of $n(f(\bar{x}(k+1))-f^{\star})$
in eq.~(D.12).
It turns out that the upper bound includes the term $\|x(k) - \bar{x}(k)\|^2$.
We then use this term together with $\|u(k) - \bar{u}(k)\|^2$, and $n(f(\bar{x}(k+1)) - f^{\star})$ as a state vector.
By analyzing interrelationships of these terms (see eqs.~(D.4)--(D.14)), we construct a novel linear system of inequalities as given by eq.~(D.25).
Please refer to the proofs of Proposition~1 and Theorems~1 and~2 for the detailed
analysis for linear convergence of Algorithm~1.
\end{remark}
}

\vspace{-0.3cm}
\section{Distributed Proportional Integral Algorithm with Finite Data Rates}\label{sec-4}
\vspace{-0.3cm}
In this section, based on the encoder--decoder scheme, we propose a distributed proportional integral algorithm with quantized communication to solve the distributed nonconvex optimization problem \LX{under the P--{$\L$} condition}, which is inspired by the proportional-integral control strategy \cite{wang2010AAC,Kia2015Aut,Gharesifard2014TAC}.
More specifically, the quantized distributed proportional integral algorithm is shown in Algorithm~2.

\floatname{algorithm}{Algorithm}
\renewcommand{\thealgorithm}{2}
\begin{algorithm}[!h]
  \caption{Quantized Distributed Proportional Integral Algorithm}
  \label{alg:A}
  \begin{algorithmic}
  \blue{
  \STATE  $\bm{\mathrm{For~each~agent}}$ $i\in\mathcal{V}.$
	\STATE $\bm{\mathrm{Initialization:}}$\\
     $x_{i}(0)\in\mathbb{R}^{m},~\sum_{j=1}^{n}u_{j}(0)=\bm{0}_{m},~\hat{x}_{j}(0)=\mathbf{0}_m$.}
  \STATE  $\bm{\mathrm{for}}$ $k\geq0:$
	\STATE $\bm{\mathrm{Update:}}$
\begin{subequations}\label{Alg}
\begin{align}
  x_{i}(k+1)&=x_{i}(k)-\xi\sum_{j=1}^{n}L_{ij}\hat{x}_{j}(k)-\varphi u_{i}(k)\notag\\
  &\quad-\sigma\nabla f_{i}(x_{i}(k)),\label{Alga}\\
  u_{i}(k+1)&=u_{i}(k)+\varphi\sum_{j=1}^{n}L_{ij}\hat{x}_{j}(k),\label{Algb}
  \end{align}
\end{subequations}
where $\sigma>0$ is the fixed step-size, $\xi$ and $\varphi$ are gain parameters.\\

\STATE $\bm{\mathrm{Send:}}$\\
The quantized output $z_{i}^{x}(k+1)$ generated by encoder~(\ref{encoder2}) to its neighbors.
%
%
\STATE $\bm{\mathrm{Receive:}}$\\
The quantized output $z_{j}^{x}(k+1)$ generated by encoder~(\ref{encoder2}) from its neighbors.
\STATE $\bm{\mathrm{Compute:}}$\\
The variable $\hat{x}_{i}(k)$ generated by decoder~(\ref{Decoder2A}).
\end{algorithmic}
\end{algorithm}

\xl{
\begin{remark}\label{rem-6}
Algorithms~1 and 2 combine the distributed gradient tracking algorithm \cite{Nedic2017SIAM,Qu2017TCNS} and the distributed proportional integral algorithm \cite{Kia2015Aut,Gharesifard2014TAC} with the quantization scheme \cite{Lei2020SIAM,Ma2021Springer}, respectively. The operational differences between Algorithm~1 and Algorithm~2 stem from the differences between the distributed gradient tracking algorithm and the distributed proportional integral algorithm. It is known from \cite{Nedic2017SIAM,Kia2015Aut,Yang2019ARC} that the distributed gradient tracking algorithm uses an auxiliary variable to track the average gradient and performs a distributed inexact gradient method, whereas the distributed proportional integral algorithm incorporates an integral feedback mechanism to correct errors caused by the distributed gradient descent method with a fixed step size.
\end{remark}}

\blue{Compared with Algorithm~\ref{alg:B}, which requires two parameters $\beta$ and $\delta$, Algorithm~\ref{alg:A} requires three parameters $\xi$, $\varphi$ and $\sigma$.}
However, in Algorithm~\ref{alg:B}, at each iteration each agent $i$ needs to communicate one additional $m$-dimensional variable besides the communication of $z_{i}^{x}(k)$ with its neighbors, which makes the quantization scheme more involved.

The following proposition provides a sufficient condition for the nonsaturation of the designed quantizer. The proof is based on the following Lyapunov candidate function:
\begin{align}
W(k)&=V(k)+n(f(\bar{x}(k))-f^{\star}),\label{Lyap-algo2}\\
V(k)&=\bm{x}^{T}(k)\bm{K}\bm{x}(k)+2\bm{x}^{T}(k)\bm{K}\bm{P}(\bm{u}(k)+\frac{\sigma}{\varphi}\bm{g}(k)) \notag \\
&\quad+(\bm{u}(k)+\frac{\sigma}{\varphi}\bm{g}(k))^{T}(\frac{\varphi+\xi}{\varphi}\bm{P})(\bm{u}(k)+\frac{\sigma}{\varphi}\bm{g}(k)), \notag
\end{align}
where $\bm{P}=P\otimes\bm{\mathrm{I}}_{m}$, $\bm{K}=K_{n}\otimes\bm{\mathrm{I}}_{m}$
with $P$ and $K_n$ are given in Lemma~\ref{lem-2},
and $\bm{g}(k)=\nabla F(\bm{x}(k))$.

\begin{proposition}\label{prop}(nonsaturation).
Suppose that \blue{Assumptions~\ref{assum-1}--\ref{assum-3}} hold. Let each agent $i\in\mathcal{V}$ run Algorithm~\ref{alg:A}, and the parameters are given as follows:
  \begin{align*}
  &\xi\in[\frac{5}{\underline{\rho}(L)}\varphi,\kappa_{1}\varphi],~\varphi\in[\sigma\kappa_{2},\sigma\kappa_{3}],\\ &\sigma\in(0,\min\{\blue{\frac{\varepsilon}{\eta_{1}},\frac{\varepsilon}{\eta_{2}}},\frac{2}{\nu},\frac{1}{4L_{f}}\}),
  \end{align*}
where the parameters $\varepsilon\in(0,\min\{\frac{\kappa_{2}}{2}-2-3L_{f}^{2}\kappa_{1}^{2}-\frac{L_{f}^{2}}{2},\kappa_{2}-1-\frac{3L_{f}^{2}+8}{\underline{\rho}(L)}\})$, $\kappa_{1}>\frac{5}{\underline{\rho}(L)}$, $\kappa_{2}>\max\{6L_{f}^{2}(\kappa_{1}+1)^{2}\kappa_{1}^{2}\rho(L),4+6L_{f}^{2}\kappa_{1}^{2}+L_{f}^{2},6L_{f}^{2}(\kappa_{1}+1)^{2},1+\frac{3L_{f}^{2}+8}{\underline{\rho}(L)}\}$ and $\kappa_{3}>\kappa_{2}$ with
\begin{align*}
\eta_{1}&=\kappa_{3}^{2}\rho(L)+\frac{2}{\underline{\rho}(L)}+2\kappa_{3}^{2}\rho(L)\\
&\quad+3\kappa_{3}^{2}L_{f}^{2}(\frac{\kappa_{1}+1}{\kappa_{2}^{2}}+\frac{3}{2}\rho(L)),\\
\eta_{2}&=4\kappa_{1}^{2}\kappa_{3}^{2}\rho^{2}(L)+2(\kappa_{3}^{2}(\kappa_{1}+1)\rho(L)+1+\kappa_{3}^{2})\\
&\quad+3\kappa_{1}^{2}L_{f}^{2}((\kappa_{1}+1)\rho(L)+\frac{3}{2}\kappa_{3}^{2}\rho^{2}(L)).
\end{align*}
Then, for any
\begin{equation}\label{K-level}
\LX{\mathcal{K}\geq\Omega},
\end{equation}
with $\mu\in(\sqrt{\epsilon_{3}},1)$ and
\begin{align*}
\Omega&=\LX{\epsilon_{1}\sqrt{\frac{\epsilon_{2}nm}{4\mu^{2}(\mu^{2}-\epsilon_{3})}}+\frac{(1+2\xi d)}{2\mu}-\frac{1}{2},}\\
\epsilon_{1}&=\max\{\xi^{2}\rho^{2}(L),\frac{\varphi^{3}\rho(L)}{\varphi+\xi},\xi\varphi\rho^{2}(L)\},\\
\epsilon_{2}&=\xi\rho(L)+2\varphi\rho(L)+4\xi^{2}\rho^{2}(L)+2(\varphi(\xi+\varphi)\rho(L)\\
&\quad+\sigma^{2}+\varphi^{2}+2\varphi),\\
\epsilon_{3}&=1-\frac{\epsilon_{4}}{\epsilon_{5}},
~\epsilon_{4}=\min\{\epsilon_{6},\epsilon_{7},\frac{\sigma}{2}\nu\},\\
\epsilon_{5}&=\max\{\frac{\xi\underline{\rho}(L)+\varphi}{\xi\underline{\rho}(L)},~1+\frac{2\xi}{\varphi}\},\\
\epsilon_{6}&=\varphi-\frac{8\sigma}{\underline{\rho}(L)}-\frac{6\sigma^{2}\varphi^{2}L_{f}^{2}(\xi+\varphi)^{2}}{\varphi^{5}}-\frac{3\sigma L_{f}^{2}}{\underline{\rho}(L)}\\
&\quad-(\varphi^{2}\rho(L)+\frac{2\sigma^{2}}{\underline{\rho}(L)}+2\varphi^{2}\rho(L)\\
&\quad+3\varphi^{2}L_{f}^{2}(\frac{\sigma^{2}(\xi+\varphi)}{\varphi^{3}}+\frac{3}{2}\rho(L))),\\
\epsilon_{7}&=\epsilon_{8}-\frac{\sigma}{2}L_{f}^{2},\\
\epsilon_{8}&=\xi\underline{\rho}(L)-\frac{9\varphi}{2}-\sigma-\frac{6\sigma^{2}\xi^{2}L_{f}^{2}(\xi+\varphi)^{2}}{\varphi^{5}}\rho(L)-\frac{3\sigma L_{f}^{2}\xi^{2}}{\varphi^{2}}\\
&\quad-(4\xi^{2}\rho^{2}(L)+2(\varphi(\xi+\varphi)\rho(L)+\sigma^{2}+\varphi^{2})\\
&\quad+3\xi^{2}L_{f}^{2}(\frac{\sigma^{2}(\xi+\varphi)}{\varphi^{3}}\rho(L)+\frac{3}{2}\rho^{2}(L))),
\end{align*}
the quantizer (\ref{Encoder2Aa}) is never saturated provided that
\begin{equation}\label{ScaF}
s(0)\geq\max\Bigg\{\frac{C_{x}+\varphi C_{u}+\sigma C_{g}}{\mathcal{K}+\frac{1}{2}},\sqrt{\frac{4\mu^{2}(\mu^{2}-\epsilon_{3})W(0)}{\epsilon_{2}nm}}\Bigg\},
\end{equation}
where $C_{x}\geq\|\bm{x}(0)\|_{\infty},~C_{u}\geq\|\bm{u}(0)\|_{\infty},~C_{g}\geq\|\bm{g}(0)\|_{\infty}$.
\end{proposition}
\textbf{Proof:} The proof is given in Appendix \ref{Appendix-A}.


\begin{theorem}\label{thm-1}
(high data rate). Suppose that Assumptions~\ref{assum-1}--\ref{assum-3} hold.
Let each agent $i\in\mathcal{V}$ run the Algorithm~\ref{alg:A} with the same $\xi$, $\varphi$, $\sigma$, $\mu$, $\mathcal{K}$ and $s(0)$ given in Proposition \ref{prop}.
Then,
\begin{equation}\label{Thm1A}
\|\bm{x}(k)-\bar{\bm{x}}(k)\|^{2}+n(f(\bar{x}(k))-f^{\star})\leq\epsilon_{9}\mu^{2k},
\end{equation}
where $\epsilon_{9}=\frac{nm\epsilon_{2}s^{2}(0)}{4\epsilon_{10}\mu^{2}(\mu^{2}-\epsilon_{3})}$, $\epsilon_{10}=\min\{\frac{\xi\underline{\rho}(L)-\varphi}{\xi\underline{\rho}(L)},1\}$.
\end{theorem}
\textbf{Proof:} The proof is given in Appendix \ref{Appendix-B}.

Similar to Theorem~\ref{thm-4}, we then have the following linear convergence result for Algorithm~\ref{alg:A}.

\begin{theorem}\label{thm-2}
(low data rate). Suppose that Assumptions~\ref{assum-1}--\ref{assum-3} hold. Let each agent $i\in\mathcal{V}$ run the Algorithm~\ref{alg:A} with the same $\xi$, $\varphi$ given in Proposition \ref{prop} and $(\mu,\sigma)\in\bar{\Pi}$, where
\begin{align*}
&\bar{\Pi}=\{(\mu,\sigma):\sigma\in(0,\min\{\blue{\frac{\varepsilon}{\eta_{1}},\frac{\varepsilon}{\eta_{2}}},\frac{2}{\nu},\frac{1}{4L_{f}}\}),\\
&\quad\quad~~\mu\in(\sqrt{\epsilon_{3}},1),~\Omega\leq \mathcal{K}\}.
\end{align*}
Then, for any $\mathcal{K}\geq1$ and $s(0)$ satisfying (\ref{ScaF}) in Proposition~\ref{prop}, \xl{$\bar{\Pi}$ is nonempty, and}
\begin{equation}\label{Thm2A}
\|\bm{x}(k)-\bar{\bm{x}}(k)\|^{2}+n(f(\bar{x}(k))-f^{\star})\leq\epsilon_{9}\mu^{2k}.
\end{equation}
\end{theorem}
\textbf{Proof:} The proof is given in Appendix \ref{Appendix-C}.

\xl{
\begin{remark}\label{rem-4}
Note that the quadratic Lyapunov used in \cite{Gharesifard2014TAC,Kia2015Aut} for convergence analysis relies on the strong convexity condition and the prefect communication.
However, such analysis cannot be used for the P--{\L} condition and quantized communication.
To tackle this problem, we design a novel Lyapunov function given in \eqref{Lyap-algo2}.
Please refer to the proofs of Proposition~2 and Theorems~3 and 4 for the detailed analysis for linear convergence of Algorithm~2.
\end{remark}

\begin{remark}
Note that \cite{zhao2022NIPS} proposed the BEER algorithm for distributed nonconvex optimization under compressed communication.
The uniform quantizer used in our proposed algorithms differs from the randomized compression operator used in \cite{zhao2022NIPS}.
Theorems~2 and~4 establish linear convergence even for 1-bit data rate under the P--{\L} condition, which is most communication efficient.
By using the randomized compression operator, \cite{zhao2022NIPS} also established the sublinear convergence for the general nonconvex case.
\end{remark}

}

\section{Numerical Examples}\label{sec-5}
\xl{In this section, we demonstrate the effectiveness of the proposed quantized distributed algorithms through two simulation studies.
In the first case, we compare the proposed algorithms with their unquantized counterparts.
In the second case, we compare the proposed algorithms with existing quantized distributed optimization algorithms.}

\begin{figure}[!hbt]
\centering
\includegraphics[scale=0.6]{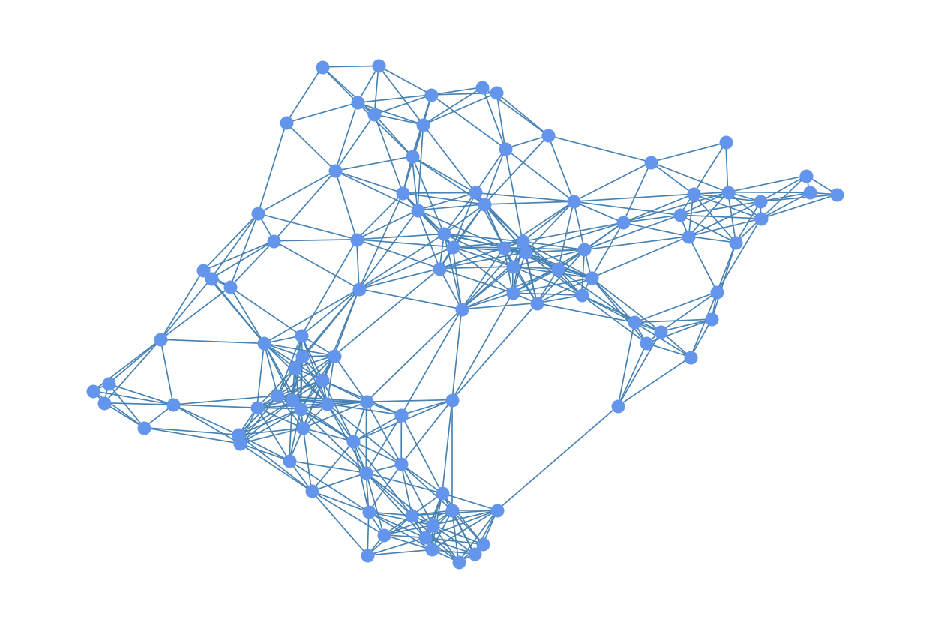}
\caption{Random connected network of 100 agents.}
\label{Graph}
\end{figure}

First, consider an undirected \blue{connected} network consisting of $100$ agents and the communication graph is randomly generated as shown in Fig. \ref{Graph}.

The local nonconvex cost functions are given by
\begin{align*}
f_i(x)&=a_{i,1}x^{2}+a_{i,2}\sin^{2}(x)+a_{i,3}\cos^{2}(x)+a_{i,4}\sin(x)\\
\quad&+a_{i,5}\sqrt{x^{4}+3}+a_{i,6}(x^{2}+2)^{1/3}+\frac{a_{i,7}x^{2}}{(\sqrt{x^{2}+1})^{1/2}}-1,
\end{align*}
where $\sum_{i=1}^{n}a_{i,1}=1$, $\sum_{i=1}^{n}a_{i,2}=4$, $\sum_{i=1}^{n}a_{i,3}=1$, $\sum_{i=1}^{n}a_{i,4}=0$, $\sum_{i=1}^{n}a_{i,5}=0$, $\sum_{i=1}^{n}a_{i,6}=0$, $\sum_{i=1}^{n}a_{i,7}=0$.
It is easy to check that Assumptions~\ref{assum-1}--\ref{assum-2} are satisfied. Moreover, it can be \LX{found} that the global cost function is \blue{$\frac{1}{100}(x^{2}+3\sin^{2}(x))$}, which satisfies Assumption~\ref{assum-3}, \blue{as shown in \cite{karimi2016linear}}.

\begin{figure}[!hbt]
\centering
\includegraphics[scale=0.61]{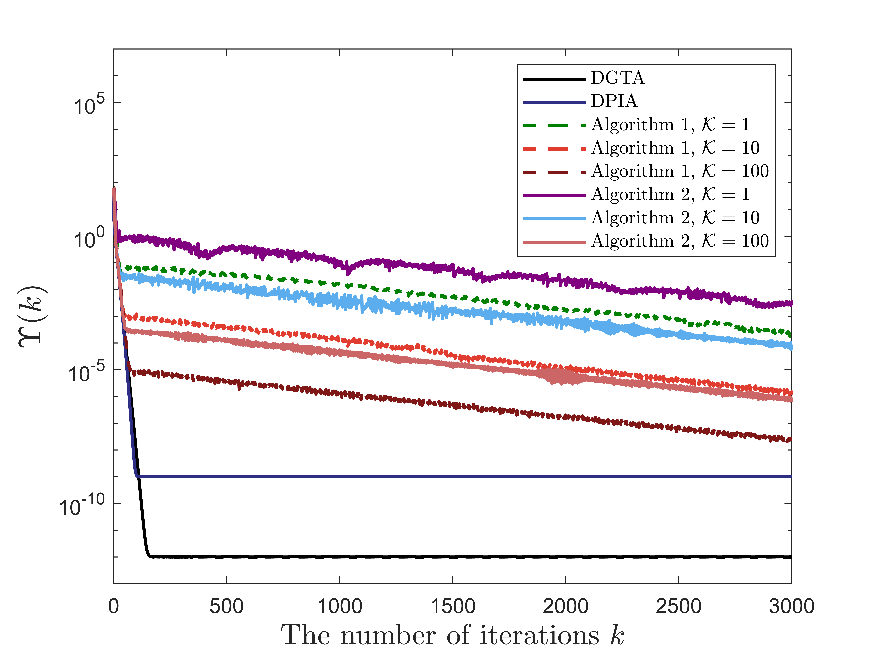}
\caption{\blue{Evolutions of $\LX{\Upsilon(k)}$ with respect to the number of iterations \red{for different algorithms}.}}
\label{Convergence}
\end{figure}

\begin{figure}[!hbt]
\centering
\includegraphics[scale=0.61]{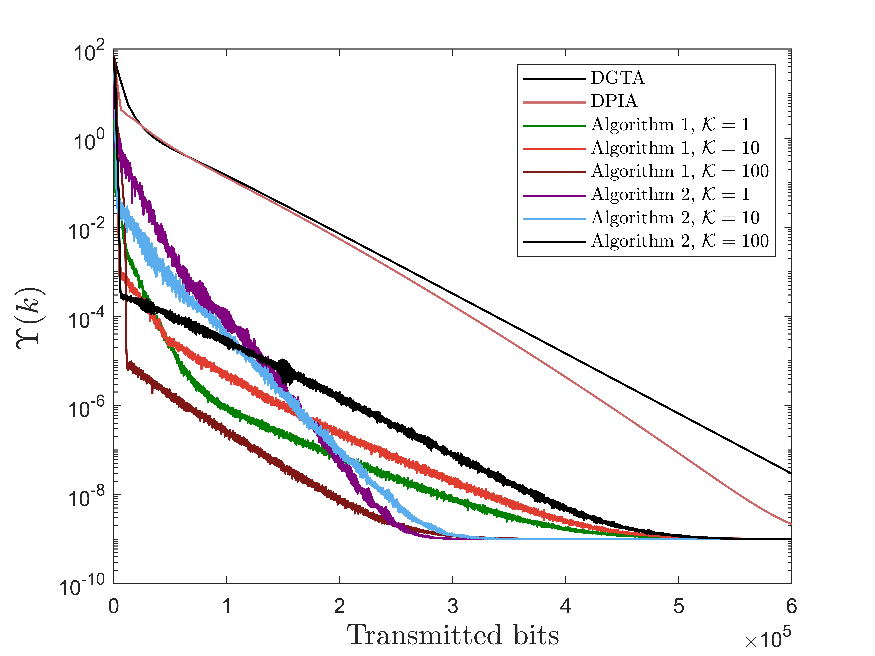}
\caption{\blue{Evolutions of $\LX{\Upsilon(k)}$ with respect to the number of transmitted bits \red{for different distributed algorithms}.}}
\label{Bit_CONVERGENCE}
\end{figure}

\begin{figure}[!hbt]
\flushright
\includegraphics[width=9.4cm,height=6.5cm]{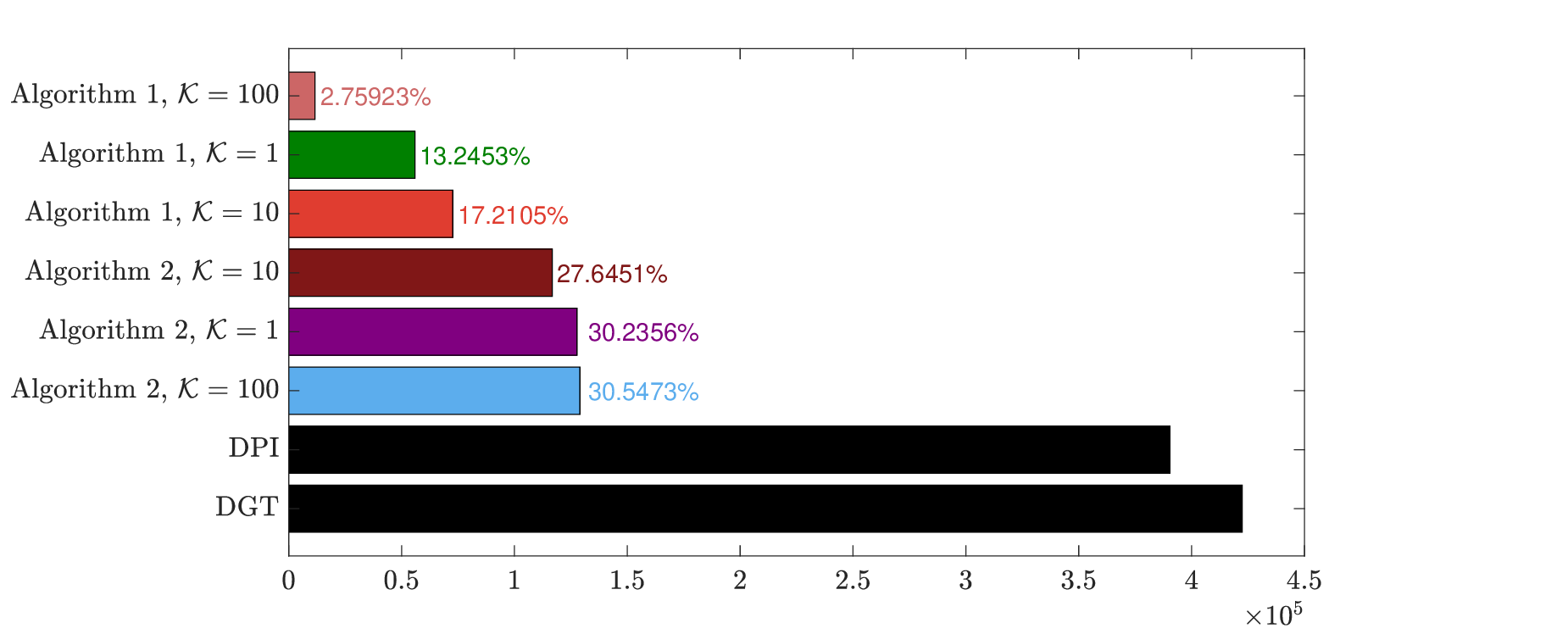}
\caption{\blue{Transmitted bits for different \red{algorithms} and quantization \red{levels} to reach $\LX{\Upsilon(k)}\leq10^{-5}$.}}
\label{Bit}
\end{figure}

\LX{Consider} the cases $\mathcal{K}=1,10,100$. \LX{Based} on the conditions (\ref{K-level}), (\ref{K-level2}), we set $s(0)=0.3198,0.0545,0.0055$, respectively.
\LX{To simplify notation, we let $\Upsilon(k)=\sum_{i=1}^{n}\|x_{i}(k)-\bar{x}(k)\|^{2}+n(f(\bar{x}(k))-f^{\star})$.}
\blue{Fig.~\ref{Convergence} and Fig.~\ref{Bit_CONVERGENCE} illustrate the \red{evolution} of $\LX{\Upsilon(k)}$ with respect to both the number of iterations $k$ and bits transmitted, \red{respectively,} for \LX{the} distributed gradient tracking algorithm (DGTA) proposed in \cite{Nedic2017SIAM}, distributed proportional integral algorithm (DPIA) proposed in \cite{wang2010AAC}, Algorithms \ref{alg:B} and \ref{alg:A}.
\LX{The} comparison of transmitted bits for the different \red{algorithms} and quantized \red{levels} to reach $\LX{\Upsilon(k)}\leq10^{-5}$ is provided in Fig.~\ref{Bit}.
The algorithm parameters used in the experiment are given in Table \ref{tbl:table2.3}.

Fig.~\ref{Convergence} clearly shows that the proposed algorithms have comparable convergence speeds as the corresponding algorithm with perfect communication, even \LX{when} the exchanged information is \LX{only} one bit.
\LX{Moreover,} larger quantization level leads to faster convergence.
This result is reasonable since a larger quantization level implies a smaller quantization error.
From Fig.~\ref{Bit_CONVERGENCE}, we can see that our proposed algorithms converge significantly faster than the DPIA and DGTA when comparing their performances based on the number of bits that each agents communicates, which shows the superiority of our proposed algorithms.}
Fig.~\ref{Bit} \LX{illustrates} that our algorithms \LX{require} only a small number of bits compared to DPIA and DGTA to reach a specific level of error.
\LX{As shown in Fig.~\ref{Bit}, for K = 100, Algorithm~\ref{alg:B} only requires $2.75923\%$ of the bits used by DGT to reach a specific level of error.}

\begin{table}[!hbt]
\caption{Parameter Settings for Our Proposed Algorithms.}
\vspace{3mm}
\centering
{
\begin{tabular}{|c|c|c|c|c|c|c|}
\hline
Algorithm & $\xi$ & $\varphi$ & $\sigma$ & $\beta$ & $\delta$ & $\mu$\\
\hline
Algorithm 1 & -- & -- & -- & 0.1 & 0.1 & 0.999\\
\hline
Algorithm 2 & 0.235 & 0.2 & 0.1 & -- & -- & 0.999\\
\hline
\end{tabular}}
\label{tbl:table2.3}
\end{table}

\xl{Next, consider an undirected connected network consisting of 10 agents and the communication graph is randomly generated.
The local cost functions associated with agents are
\begin{align*}
f_i(x)=\alpha_{i}\frac{(x^{1}-\sin(x^{2}))^{2}}{2},
\end{align*}
where $i\in\{1,\cdots,10\}$, $\alpha_{1}=\alpha_{4}=\alpha_{7}=0.1$, $\alpha_{2}=\alpha_{6}=\alpha_{8}=\alpha_{9}=0.05$, and $\alpha_{3}=\alpha_{10}=0.15$, $x=[x^{1},x^{2}]^{T}$.
This function is commonly used in deep learning applications and satisfies the P--{\L} condition \cite{apidopoulos2022JGO}.
It is easy to check that Assumptions~\ref{assum-1}--\ref{assum-3} are satisfied.
Fig.~\ref{CompareSevAlg1} illustrates the evolution of $\Upsilon(k)$ with respect to the number of iterations $k$, for different quantized distributed optimization algorithms proposed in \cite{Yi2014TCNS,Kajiyama2021TAC}, and Algorithms~\ref{alg:B} and \ref{alg:A}.
Choose $\mathcal{K}=300$ and $s(0)=1$ such that the conditions (\ref{K-level2}) and (\ref{K-level}) are satisfied.
The parameters of the various algorithms used in the experiments are provided in Table \ref{tbl}.

\begin{table*}[!hbt]
\centering
\xl{
\caption{Parameter Settings for Different Quantized Distributed Algorithms.}
\centering
{
\begin{tabular}{|c|c|c|c|c|c|c|c|c|}
\hline
Algorithm & $\xi$ & $\varphi$ & $\sigma$ & $\beta$ & $\delta$ & $h$ & $\eta$ & $\mu$ \\
\hline
Algorithm 1 & -- & -- & -- & 0.2 & 0.2 & -- & -- & 0.99 \\
\hline
Algorithm 2 & 0.01 & 0.04 & 0.1 & -- & -- & -- & -- & 0.99 \\
\hline
\cite{Yi2014TCNS} & -- & -- & -- & -- & -- & 0.5 & -- & 0.99  \\
\hline
\cite{Kajiyama2021TAC} & -- & -- & -- & -- & -- & -- & 0.01 & 0.99\\
\hline
\end{tabular}}
\label{tbl}}
\end{table*}

Fig.~\ref{CompareSevAlg1} shows that Algorithms~1 and 2 are faster than the quantized subgradient descent algorithm in \cite{Yi2014TCNS} and the quantized consensus-based algorithm in \cite{Kajiyama2021TAC}.

\begin{figure}[!hbt]
\centering\xl{
\includegraphics[scale=0.61]{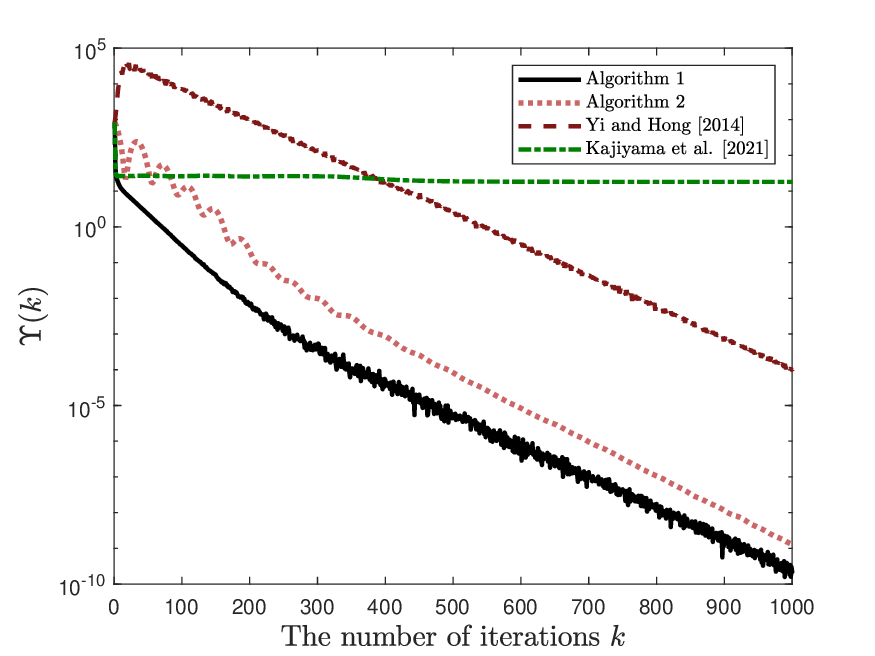}
\caption{Evolutions of $\Upsilon(k)$ with respect to the number of iterations for different quantized distributed algorithms.}
\label{CompareSevAlg1}}
\end{figure}
}

\section{Conclusions}\label{sec-6}
\LX{
In this paper, we introduced an encoder--decoder scheme to reduce the number of transmitted bits. By integrating it with distributed gradient tracking \LX{and} distributed proportional integral algorithms, respectively, we then proposed two quantized distributed algorithms for solving nonconvex optimization over an undirected connected network}.
For the case where local cost functions are smooth and the global cost function satisfies the P--{\L} condition, we showed that the proposed algorithms linearly converge to a global optimal point provided that the quantization level is larger than a certain threshold.  
\blue{We also} showed that, with appropriate algorithm parameters, the proposed algorithms with a low data rate, \LX{even bit rate one}, are sufficient to ensure linear convergence.
One future direction is to consider directed graphs.

\bibliography{Ref}

\appendix

\section{Useful Lemma}\label{Appendix-Lem1}
The following lemma \LX{is} used in this paper.
\begin{lemma}\label{lem-2}
(\citet[Lemma~3]{Yi2022Arxiv})
Let $L$ be the Laplacian matrix of an undirected and connected graph $\mathcal{G}$ with $n$ agents and $K_{n}=\bm{\mathrm{I}}_{n}-\frac{1}{n}\bm{1}_{n}\bm{1}_{n}^{T}$. Then $L$ and $K_{n}$ are positive semi-definite, $L\leq\rho(L)\bm{\mathrm{I}}_{n}$, $\rho(K_{n})=1$,
  \begin{subequations}
  \begin{align}
  &K_{n}L=LK_{n}=L,\label{Lemma1Aa}\\
  &0\leq\underline{\rho}(L)K_{n}\leq L\leq\rho(L)K_{n}.\label{Lemma1Ab}
  \end{align}
  \end{subequations}
Moreover, there exists an orthogonal matrix $[\begin{array}{cc} r & R \end{array}]\in\mathbb{R}^{n\times n}$ with $r=\frac{1}{\sqrt{n}}\bm{1}_{n}$ and $R\in\mathbb{R}^{n\times(n-1)}$ such that
  \begin{subequations}
  \begin{align}
  &PL=LP=K_{n},\label{Lemma1Bb}\\
  &\frac{1}{\rho(L)}I_{n}\leq P\leq\frac{1}{\underline{\rho}(L)}I_{n},\label{Lemma1Bc}
  \end{align}
  \end{subequations}
where $\Lambda_{1}=\rm{diag}([\lambda_{2},\dots,\lambda_{n}])$ with $0<\lambda_{2}\leq\dots\leq\lambda_{n}$ beging the nonzero eigenvalues of the Laplacian matrix $L$, and
  \begin{align*}
  &P=\left[\begin{array}{cc}
           r & R
         \end{array}\right]
    \left[\begin{array}{cc}
      \lambda^{-1}_{n} & 0 \\
      0 & \Lambda^{-1}_{1} \\
    \end{array}\right]
  \left[\begin{array}{c}
    r^{T} \\
    R^{T}
  \end{array}\right].
  \end{align*}
\end{lemma}

\section{Proof of Lemma \ref{lem-1}}\label{Appendix-Lem1}
First, from $\beta\in(0,\LX{\frac{\sqrt{2}}{2\rho(L)}})$, it can be guaranteed that $\varrho=\rho(\bm{\mathrm{I}}_{nm}-\beta\bm{L}-\bm{H})\in(0,1)$. Then, from $\delta\in(0,\frac{\sqrt{c_{1}\Theta_{1}}}{2})$, $\LX c_{2}\in(0,1-\varrho^{2})$, we have
\begin{equation}\label{IInequality1}
2\delta^{2}(1+\frac{1}{\sigma_{1}})<\frac{(1-\varrho^{2}-\LX c_{2})\Theta_{1}}{2},
\end{equation}
where the inequality holds due to $1+\frac{1}{\sigma_{1}}=\frac{1+\varrho^{2}}{1-\varrho^{2}}$. Since $(1+\sigma_{1})\varrho^{2}=\frac{1+\varrho^{2}}{2}$, it can be derived that
\begin{equation}\label{IInequality7}
1-\frac{\LX c_{2}}{2}-(1+\sigma_{1})\varrho^{2}=\frac{1-\varrho^{2}-\LX c_{2}}{2}.
\end{equation}
Then, based on (\ref{IInequality1}) and (\ref{IInequality7}), it can be guaranteed that
\begin{equation}\label{Inequality1}
\chi_{1}\Theta_{1}+\chi_{2}<(1-\frac{\LX c_{2}}{2})\Theta_{1}.
\end{equation}
Second, from $\beta\in(0,\frac{\sqrt{2}}{2\rho(L)})$ and $\delta\in(0,\frac{1}{4L_{f}})$, we have
\begin{equation}\label{IInequality2}
\beta^{2}\rho^{2}(L)<\frac{1}{2},~\frac{2\delta^{2}L_{f}^{2}}{1-2\delta L_{f}}<\frac{1}{4}.
\end{equation}
Then, from (\ref{IInequality2}) and $\Theta_{1}\in(0,\frac{c_{1}}{24L_{f}^{2}}]$, it can be derived that
\begin{equation}\label{IInequality3}
(1+\frac{1}{\sigma_{1}})8L_{f}^{2}(\beta^{2}\rho^{2}(L)+\frac{2\delta^{2}L_{f}^{2}}{1-2\delta L_{f}})\Theta_{1}<\frac{1-\varrho^{2}-\LX c_{2}}{4}.
\end{equation}
From $\delta\in(0,\min\{\frac{2}{\nu},\frac{1}{8+2L_{f}}\})$, it can be calculated that
\begin{equation}\label{IInequality4}
\frac{2\delta(2-\delta\nu)}{1-2\delta L_{f}}<\frac{4\delta}{1-2\delta L_{f}}<\frac{1}{2}.
\end{equation}
Based on (\ref{IInequality4}), we have
\begin{equation}\label{IInequality5}
(1+\frac{1}{\sigma_{1}})\frac{\LX{16L_{f}^{2}\delta(2-\delta\nu)}}{1-2\delta L_{f}}\Theta_{2}<\frac{1-\varrho^{2}-\LX c_{2}}{8}.
\end{equation}
Then, from $\delta\in(0,\frac{\sqrt{c_{1}}}{8L_{f}})$, we have
\begin{equation}\label{IInequality6}
(1+\frac{1}{\sigma_{1}})8L_{f}^{2}\delta^{2}<\frac{1-\varrho^{2}-\LX c_{2}}{8}.
\end{equation}
Based on (\ref{IInequality7}), (\ref{IInequality3}), (\ref{IInequality5}) and (\ref{IInequality6}), it can be guaranteed that
\begin{equation}\label{Inequality7}
\chi_{3}\Theta_{1}+\chi_{4}+\chi_{5}\Theta_{2}<1-\frac{\LX c_{2}}{2}.
\end{equation}
Third, from $\Theta_{1}\in(0,\frac{\nu\Theta_{2}}{2L_{f}^{2}}]$, we know that
\begin{equation}\label{Inequality8}
\frac{\delta}{2}L_{f}^{2}\Theta_{1}\leq\frac{\delta}{4}\nu\Theta_{2}.
\end{equation}
Based on (\ref{Inequality8}) and $\delta\in(0,\frac{2}{\nu})$, it can be guaranteed that
\begin{equation}\label{Inequality9}
\chi_{6}\Theta_{1}+\chi_{7}\Theta_{2}\leq(1-\frac{\delta}{4}\nu)\Theta_{2}.
\end{equation}
From (\ref{Inequality1}), (\ref{Inequality7}) and (\ref{Inequality9}), we know that (\ref{Phi}) holds. This completes the proof.

\section{Proof of Lemma \ref{lem-3}}\label{Appendix-Lem3}
For the matrix $\Phi$ in (\ref{PPhi}), we have
  \begin{align}\label{Xi}
  (\bm{\mathrm{I}}_{3}+\Phi)^{2}
  =\left[
     \begin{array}{ccc}
       \Xi_{1} & \Xi_{2} & \Xi_{3} \\
       \Xi_{4} & \Xi_{5} & \Xi_{6} \\
       \Xi_{7} & \Xi_{8} & \Xi_{9} \\
     \end{array}
   \right],
  \end{align}
where
  \begin{align*}
  \Xi_{1}&=(1+\chi_{1})^{2}+\chi_{2}\chi_{3},~\Xi_{2}=(1+\chi_{1})\chi_{2}+\chi_{2}(1+\chi_{4}),\\
  \Xi_{3}&=\chi_{2}\chi_{5},~\Xi_{4}=\chi_{3}(1+\chi_{1})+(1+\chi_{4})\chi_{3}+\chi_{5}\chi_{6},\\
  \Xi_{5}&=\chi_{3}\chi_{2}+(1+\chi_{4})^{2},~\Xi_{6}=(1+\chi_{4})\chi_{5}+\chi_{5}(1+\chi_{7}),\\
  \Xi_{7}&=\chi_{6}(1+\chi_{1})+(1+\chi_{7})\chi_{6},~\Xi_{8}=\chi_{6}\chi_{2},\\
  \Xi_{9}&=(1+\chi_{7})^{2}.
  \end{align*}
Based on (\ref{Xi}) and Lemma~\ref{lem-1}, it can be found that \LX{$(\bm{\mathrm{I}}_{3}+\Phi)^{2}\succ0$}. This together with \citet[Lemma 8.4.1]{Horn2012Matrix} implies $\Phi$ is irreducible. Since $\Phi$ is irreducible and nonnegative, by \citet[Lemma 8.4.4]{Horn2012Matrix}, $\rho(\Phi)>0$ and there exists a positive vector $\zeta$ such that $\Phi\zeta=\rho(\Phi)\zeta$. Then, by \citet[Corollary 8.1.33]{Horn2012Matrix}, we obtain
\begin{equation}\label{PhiIneq}
\LX{\Phi^{k}\preceq\rho^{k}(\Phi)\frac{\max_{1\leq i\leq 3}\zeta_{i}}{\min_{1\leq i\leq 3}\zeta_{i}}\mathbf{1}_{3}\mathbf{1}_{3}^{T}.}
\end{equation}
Hence, (\ref{PhiIneq}) yields (\ref{hatrho}). This proof is complete.

\section{Proof of Proposition \ref{prop2}}\label{Appendix-D}

Define $e_{j}^{x}(k)=x_{j}(k)-\hat{x}_{j}(k)$, $e_{j}^{u}(k)=u_{j}(k)-\hat{u}_{j}(k)$. Since the initial value $b_{j}^{x}(0)=\hat{x}_{j}(0)$, $b_{j}^{u}(0)=\hat{u}_{j}(0)$, we obtain that $b_{j}^{x}(k)=\hat{x}_{j}(k)$ and $b_{j}^{u}(k)=\hat{u}_{j}(k)$. Hence, we have $e_{j}^{x}(k)=x_{j}(k)-b_{j}^{x}(k)$, $e_{j}^{u}(k)=u_{j}(k)-b_{j}^{u}(k)$.
Denote $\bm{e}^{x}(k)=[(e^{x}_{1}(k))^T,\dots,(e^{x}_{n}(k))^T]$, $\bm{e}^{u}(k)=[(e^{u}_{1}(k))^T,\dots,(e^{u}_{n}(k))^T]$.
\LX{Then, the update rule (\ref{Alg2}) can be rewritten in a compact form:}
  \begin{subequations}
  \begin{align}
  \bm{x}(k+1)&=\bm{x}(k)-\beta\bm{L}\bm{x}(k)+\beta\bm{L}\bm{e}^{x}(k)-\delta\bm{u}(k),\label{CptAlg2A}\\
  \quad\quad~\bm{u}(k+1)&=\bm{u}(k)-\beta\bm{L}\bm{u}(k)+\beta\bm{L}\bm{e}^{u}(k)\notag\\
  &\quad+\nabla F(\bm{x}(k+1))-\nabla F(\bm{x}(k)),\notag\\
  &\quad~\forall\bm{x}(0)\in\mathbb{R}^{nm},~\bm{u}(0)=\nabla F(\bm{x}(0)).\label{CptAlg2B}
  \end{align}
  \end{subequations}

Denote $\bm{g}(k)=\nabla F(\bm{x}(k))$, $\bar{\bm{g}}(k)=\bm{H}\bm{g}(k)$, $\bm{g}^{0}(k)=\nabla F(\bar{\bm{x}}(k))$, $\bar{\bm{g}}^{0}(k)=\bm{H}\bm{g}^{0}(k)=\bm{1}_{n}\otimes\nabla f(\bar{x}(k))$, $\bar{\bm{u}}(k)=\bm{1}_{n}\otimes\bar{u}(k),$
and $\bar{u}(k)=\frac{1}{n}(\bm{1}_{n}^{T}\otimes\bm{\mathrm{I}}_{m})\bm{u}(k),$
then from (\ref{CptAlg2B}) and $\bm{u}(0)=\nabla F(\bm{x}(0))$, we have
\begin{align}\label{barU2}
\bar{\bm{u}}(k+1)&=H(\bm{u}(k)-\beta\bm{L}\bm{u}(k)+\beta\bm{L}\bm{e}^{u}(k)\notag\\
  &\quad+\bm{g}(k+1)-\bm{g}(k))\notag\\
&=\bar{\bm{g}}(k+1).
\end{align}

From (\ref{CptAlg2A}), it can be calculated that
  \begin{align}\label{barX2}
  \bar{\bm{x}}(k+1)&=H(\bm{x}(k)-\beta\bm{L}\bm{x}(k)+\beta\bm{L}\bm{e}^{x}(k)-\delta\bm{u}(k))\notag\\
  &=\bar{\bm{x}}(k)-\delta\bar{\bm{u}}(k).
  \end{align}

Based on (\ref{CptAlg2A}) and (\ref{barX2}), we have
  \begin{align}\label{ErrorX2}
  &\|\bm{x}(k+1)-\bar{\bm{x}}(k+1)\|^{2}\notag\\
  &=\|\bm{x}(k)-\beta\bm{L}\bm{x}(k)+\beta\bm{L}\bm{e}^{x}(k)-\delta\bm{u}(k)\notag\\
  &\quad-(\bar{\bm{x}}(k)-\delta\bar{\bm{u}}(k))\|^{2}\notag\\
  &\LX{=\|(\bm{\mathrm{I}}_{nm}-\beta\bm{L}-\bm{H})(\bm{x}(k)-\bar{\bm{x}}(k))+\beta\bm{L}\bm{e}^{x}(k)}\notag\\
  &\quad\LX{-\delta(\bm{u}(k)-\bar{\bm{u}}(k))\|^{2}}\notag\\
  &\leq(1+\sigma_{1})\varrho^{2}\|\bm{x}(k)-\bar{\bm{x}}(k)\|^{2}\notag\\
  &\quad+(1+\frac{1}{\sigma_{1}})\|\beta\bm{L}\bm{e}^{x}(k)-\delta(\bm{u}(k)-\bar{\bm{u}}(k))\|^{2}\notag\\
  &\leq(1+\frac{1}{\sigma_{1}})(2\beta^{2}\rho^{2}(L)\|\bm{e}^{x}(k)\|^{2}+2\delta^{2}\|\bm{u}(k)-\bar{\bm{u}}(k)\|^{2})\notag\\
  &\quad+(1+\sigma_{1})\varrho^{2}\|\bm{x}(k)-\bar{\bm{x}}(k)\|^{2},
  \end{align}
where \LX{the second equality holds due to $\bm{L}\bar{\bm{x}}(k)=\bm{L}\bm{H}\bm{x}(k)=\bm{0}_{nm}$ and $\bm{H}\bm{x}(k)=\bm{H}\bar{\bm{x}}(k)$}; the first inequality holds due to the Cauchy--Schwarz inequality; the last inequality holds due to (\ref{Lemma1Ab}), $\rho(\bm{K})\leq1$ and the Cauchy--Schwarz inequality.

Based on (\ref{CptAlg2B}), we have
  \begin{align}\label{ErrorU2}
  &\|\bm{u}(k+1)-\bar{\bm{u}}(k+1)\|^{2}\notag\\
  &=\|\bm{u}(k)-\beta\bm{L}\bm{u}(k)+\beta\bm{L}\bm{e}^{u}(k)+\bm{g}(k+1)-\bm{g}(k)\notag\\
  &\quad-\bar{\bm{u}}(k+1)\|^{2}\notag\\
  &=\LX{\|(\bm{\mathrm{I}}_{nm}-\beta\bm{L}-\bm{H})(\bm{u}(k)-\bar{\bm{u}}(k))+\beta\bm{L}\bm{e}^{u}(k)}\notag\\
  &\quad\LX{+(\bm{\mathrm{I}}_{nm}-\bm{H})(\bm{g}(k+1)-\bm{g}(k))}\notag\\
  &\leq(1+\frac{1}{\sigma_{1}})\|\beta\bm{L}\bm{e}^{u}(k)+(\bm{\mathrm{I}}_{nm}-\bm{H})(\bm{g}(k+1)-\bm{g}(k))\|^{2}\notag\\
  &\quad+(1+\sigma_{1})\varrho^{2}\|\bm{u}(k)-\bar{\bm{u}}(k)\|^{2}\notag\\
  &\leq(1+\frac{1}{\sigma_{1}})(2\beta^{2}\rho^{2}(L)\|\bm{e}^{u}(k)\|^{2}+2\|(\bm{\mathrm{I}}_{nm}-\bm{H})\notag\\
  &\quad\times(\bm{g}(k+1)-\bm{g}(k))\|^{2})+(1+\sigma_{1})\varrho^{2}\|\bm{u}(k)-\bar{\bm{u}}(k)\|^{2}\notag\\
  &\leq(1+\frac{1}{\sigma_{1}})(2\beta^{2}\rho^{2}(L)\|\bm{e}^{u}(k)\|^{2}+2\|\bm{g}(k+1)-\bm{g}(k)\|^{2})\notag\\
  &\quad+(1+\sigma_{1})\varrho^{2}\|\bm{u}(k)-\bar{\bm{u}}(k)\|^{2},
  \end{align}
where \LX{the second equality holds due to (\ref{barU2}), $\bm{L}\bar{\bm{u}}(k)=\bm{L}\bm{H}\bm{u}(k)=\bm{0}_{nm}$ and $\bm{H}\bm{u}(k)=\bm{H}\bar{\bm{u}}(k)$}; the first inequality holds due to the Cauchy--Schwarz inequality; the second inequality holds due to (\ref{Lemma1Ab}), $\rho(\bm{K})\leq1$ and the Cauchy--Schwarz inequality; and the last inequality holds due to $\rho(\bm{\mathrm{I}}_{nm}-\bm{H})=1$.

Based on Assumption \ref{assum-2}, we have
\blue{
  \begin{align}\label{ErrorgReg}
  &\|\bm{g}(k+1)-\bm{g}(k)\|\notag\\
  &\leq \LX{L_{f}\|\bm{x}(k+1)-\bm{x}(k)\|}\notag\\
  &= L_{f}\|\bm{x}(k)-\beta\bm{L}\bm{x}(k)+\beta\bm{L}\bm{e}^{x}(k)-\delta\bm{u}(k)-\bm{x}(k)\|\notag\\
  &=L_{f}\|-\beta\bm{L}(\bm{x}(k)-\bar{\bm{x}}(k))+\beta\bm{L}\bm{e}^{x}(k)\notag\\
  &\quad-\delta(\bm{u}(k)-\bar{\bm{u}}(k))-\delta\bar{\bm{g}}(k)\|\notag\\
  &\leq L_{f}[\beta\rho(L)\|\bm{x}(k)-\bar{\bm{x}}(k)\|+\beta\rho(L)\|\bm{e}^{x}(k)\|\notag\\
  &\quad+\delta\|\bm{u}(k)-\bar{\bm{u}}(k)\|+\delta\|\bar{\bm{g}}(k)\|],
  \end{align}
where \LX{the second equality holds due to $\bm{L}\bar{\bm{x}}(k)=\bm{L}\bm{H}\bm{x}(k)=\bm{0}_{nm}$ and (\ref{barU2})}; and the last inequality holds due to (\ref{Lemma1Ab}) and $\rho(\bm{K})\leq1$.

Based on (\ref{ErrorgReg}) and the Cauchy--Schwarz inequality, we have

  \begin{align}\label{Errorg2}
  &\|\bm{g}(k+1)-\bm{g}(k)\|^{2}\notag\\
  &\leq L_{f}^{2}\|\bm{x}(k+1)-\bm{x}(k)\|^{2}\notag\\
  &= L_{f}^{2}[4\beta^{2}\rho^{2}(L)\|\bm{x}(k)-\bar{\bm{x}}(k)\|^{2}+4\beta^{2}\rho^{2}(L)\|\bm{e}^{x}(k)\|^{2}\notag\\
  &\quad+4\delta^{2}\|\bm{u}(k)-\bar{\bm{u}}(k)\|^{2}+4\delta^{2}\|\bar{\bm{g}}(k)\|^{2}].
  \end{align}}

From Assumption \ref{assum-3}, we have
  \begin{equation}\label{barg}
  \|\bar{\bm{g}}^{0}(k)\|^{2}=n\|\nabla f(\bar{x}(k))\|^{2}\geq2\nu n(f(\bar{x}(k))-f^{\star}).
  \end{equation}

From Assumption \ref{assum-2}, it can be calculated that
  \begin{equation}\label{g0g}
  \|\bm{g}^{0}(k)-\bm{g}(k)\|^{2}\leq L_{f}^{2}\|\bar{\bm{x}}(k)-\bm{x}(k)\|^{2}.
  \end{equation}

Then, from (\ref{g0g}) and $\rho(\bm{H})=1$, we have
  \begin{align}\label{barg0g}
  \|\bar{\bm{g}}^{0}(k)-\bar{\bm{g}}(k)\|^{2}&=\|\bm{H}(\bm{g}^{0}(k)-\bm{g}(k))\|^{2}\notag\\
  &\leq\|\bm{g}^{0}(k)-\bm{g}(k)\|^{2}\notag\\
  &\leq L_{f}^{2}\|\bar{\bm{x}}(k)-\bm{x}(k)\|^{2}.
  \end{align}

From \citet[Lemma 1.2.3]{nesterov2018lectures}, we know that (\ref{assum2}) implies
\begin{equation}\label{assum2b}
|f(y)-f(x)-(y-x)^{T}\nabla f(x)|\leq\frac{L_{f}}{2}\|y-x\|^{2},~\forall x,y\in\mathbb{R}^{m}.
\end{equation}

From the P--{\L} condition, we know that
\begin{align}\label{OptimalF2}
&n(f(\bar{x}(k+1))-f^{\star})\notag\\
&=\LX{ n(f(\bar{x}(k))- f^{\star}+ f(\bar{x}(k+1))- f(\bar{x}(k)))}\notag\\
&\leq \LX{n(f(\bar{x}(k))- f^{\star})-\delta\bar{\bm{g}}^{T}(k)\bm{g}^{0}(k)+\frac{\delta^{2}L_{f}}{2}\|\bar{\bm{g}}(k)\|^{2}}\notag\\
&= \LX{n(f(\bar{x}(k))- f^{\star})-\delta\bar{\bm{g}}^{T}(k)\bar{\bm{g}}^{0}(k)+\frac{\delta^{2}L_{f}}{2}\|\bar{\bm{g}}(k)\|^{2}}\notag\\
&=n(f(\bar{x}(k))-f^{\star})-\frac{\delta}{2}\bar{\bm{g}}^{T}(k)(\bar{\bm{g}}(k)+\bar{\bm{g}}^{0}(k)-\bar{\bm{g}}(k))\notag\\
&\quad-\frac{\delta}{2}(\bar{\bm{g}}(k)+\bar{\bm{g}}^{0}(k)-\bar{\bm{g}}^{0}(k))^{T}\bar{\bm{g}}^{0}(k)+\frac{\delta^{2}L_{f}}{2}\|\bar{\bm{g}}(k)\|^{2}\notag\\
&\leq n(f(\bar{x}(k))-f^{\star})-\frac{\delta}{4}\|\bar{\bm{g}}(k)\|^{2}+\frac{\delta}{4}\|\bar{\bm{g}}^{0}(k)-\bar{\bm{g}}(k)\|^{2}\notag\\
&\quad-\frac{\delta}{4}\|\bar{\bm{g}}^{0}(k)\|^{2}+\frac{\delta}{4}\|\bar{\bm{g}}^{0}(k)-\bar{\bm{g}}(k)\|^{2}+\frac{\delta^{2} L_{f}}{2}\|\bar{\bm{g}}(k)\|^{2}\notag\\
&\leq n(f(\bar{x}(k))-f^{\star})-\frac{\delta}{4}(1-2 \delta L_{f})\|\bar{\bm{g}}(k)\|^{2}\notag\\
&\quad+\frac{\delta}{2}\|\bar{\bm{g}}^{0}(k)-\bar{\bm{g}}(k)\|^{2}-\frac{\delta}{2}\nu n(f(\bar{x}(k))-f^{\star})\notag\\
&\leq n(f(\bar{x}(k))-f^{\star})-\frac{\delta}{4}(1-2 \delta L_{f})\|\bar{\bm{g}}(k)\|^{2}\notag\\
&\quad+\frac{\delta}{2}L_{f}^{2}\|\bm{x}(k)-\bar{\bm{x}}(k)\|^{2}-\frac{\delta}{2}\nu n(f(\bar{x}(k))-f^{\star}),
\end{align}
%
where the first inequality holds since that \LX{$ f$} is smooth and (\ref{assum2b}); the \LX{second} equality holds due to $\bar{\bm{g}}^{T}(k)\bm{g}^{0}(k)=\bm{g}^{T}(k)\bm{H}\bm{g}^{0}(k)=\bm{g}^{T}(k)\bm{H}\bm{H}\bm{g}^{0}(k)=\bar{\bm{g}}^{T}(k)\bar{\bm{g}}^{0}(k)$; the second inequality holds due to the
Cauchy--Schwarz inequality;
the third inequality holds due to (\ref{barg});
and the last inequality holds due to (\ref{barg0g}).


\LX{From (\ref{OptimalF2}), we know that}
\begin{align}\label{OptimalF2c}
\|\bar{\bm{g}}(k)\|^{2}&\leq\frac{2(2-\delta\nu)}{\delta(1-2 \delta L_{f})}n(f(\bar{x}(k))-f^{\star})\notag\\
&\quad+\frac{2L_{f}^{2}}{1-2\delta L_{f}}\|\bm{x}(k)-\bar{\bm{x}}(k)\|^{2}.
\end{align}

Based on (\ref{ErrorU2}), (\ref{Errorg2}) and (\ref{OptimalF2c}), we have
\begin{align}\label{ErrorU21}
  &\|\bm{u}(k+1)-\bar{\bm{u}}(k+1)\|^{2}\notag\\
  &\leq(1+\frac{1}{\sigma_{1}})(2\beta^{2}\rho^{2}(L)\|\bm{e}^{u}(k)\|^{2}+8L_{f}^{2}\notag\\
  &\quad\times(\beta^{2}\rho^{2}(L)\|\bm{x}(k)-\bar{\bm{x}}(k)\|^{2}+\beta^{2}\rho^{2}(L)\|\bm{e}^{x}(k)\|^{2}\notag\\
  &\quad+\delta^{2}\|\bm{u}(k)-\bar{\bm{u}}(k)\|^{2}+\frac{2\delta(2-\delta\nu)}{1-2 \delta L_{f}}n(f(\bar{x}(k))-f^{\star})\notag\\
  &\quad+\frac{2\delta^{2}L_{f}^{2}}{1-2\delta L_{f}}\|\bm{x}(k)-\bar{\bm{x}}(k)\|^{2}))\notag\\
  &\quad+(1+\sigma_{1})\varrho^{2}\|\bm{u}(k)-\bar{\bm{u}}(k)\|^{2}\notag\\
  &=(1+\frac{1}{\sigma_{1}})2\beta^{2}\rho^{2}(L)(\|\bm{e}^{u}(k)\|^{2}+4L_{f}^{2}\|\bm{e}^{x}(k)\|^{2})\notag\\
  &\quad+(1+\frac{1}{\sigma_{1}})8L_{f}^{2}(\beta^{2}\rho^{2}(L)+\frac{2\delta^{2}L_{f}^{2}}{1-2\delta L_{f}})\|\bm{x}(k)-\bar{\bm{x}}(k)\|^{2}\notag\\
  &\quad+(1+\frac{1}{\sigma_{1}})\LX{\frac{16L_{f}^{2}\delta(2-\delta\nu)}{1-2 \delta L_{f}}}n(f(\bar{x}(k))-f^{\star})\notag\\
  &\quad+((1+\sigma_{1})\varrho^{2}+(1+\frac{1}{\sigma_{1}})8L_{f}^{2}\delta^{2})\|\bm{u}(k)-\bar{\bm{u}}(k)\|^{2}.
  \end{align}

Denote \LX{$\bm{b}^{x}(k)=[(b^{x}_{1}(k))^T,\dots,(b^{x}_{n}(k))^T]^T$, $\bm{b}^{u}(k)=[(b^{u}_{1}(k))^T,\dots,(b^{u}_{n}(k))^T]^T$}.
Based on (\ref{Encoder2Ba}), we know that
  \begin{subequations}\label{Encoder2A}
  \begin{align}
  \bm{b}^{x}(k)&=s(k-1)Q\left[\frac{1}{s(k-1)}(\bm{x}(k)-\bm{b}^{x}(k-1))\right]\notag\\
  &\quad+\bm{b}^{x}(k-1),\label{Estimation2a1}\\
  \bm{b}^{u}(k)&=s(k-1)Q\left[\frac{1}{s(k-1)}(\bm{u}(k)-\bm{b}^{u}(k-1))\right]\notag\\
  &\quad+\bm{b}^{u}(k-1).\label{Estimation2a2}
  \end{align}
  \end{subequations}

By subtracting $\bm{b}^{x}(k)$ from both sides of (\ref{CptAlg2A}), it can be obtained that
\begin{align}\label{XB2}
&\bm{x}(k+1)-\bm{b}^{x}(k)\notag\\
&=\bm{x}(k)-\beta\bm{L}\bm{x}(k)+\beta\bm{L}\bm{e}^{x}(k)-\delta\bm{u}(k)-\bm{b}^{x}(k)\notag\\
&=(\bm{\mathrm{I}}_{nm}+\beta\bm{L})\bm{e}^{x}(k)-\beta\bm{L}\bm{x}(k)-\delta\bm{u}(k)\notag\\
&=s(k)((\bm{\mathrm{I}}_{nm}+\beta\bm{L})\LX{\frac{\bm{e}^{x}(k)}{s(k)}}-\frac{1}{s(k)}(\beta\bm{L}\bm{x}(k)+\delta\bm{u}(k)))\notag\\
&=s(k)\theta^{x}(k),
\end{align}

where
$$
\theta^{x}(k)=(\bm{\mathrm{I}}_{nm}+\beta\bm{L})\LX{\frac{\bm{e}^{x}(k)}{s(k)}}-\frac{1}{s(k)}(\beta\bm{L}\bm{x}(k)+\delta\bm{u}(k)).
$$
Based on (\ref{Estimation2a1}) and (\ref{XB2}), we have
\begin{align}\label{Error2}
&\bm{e}^{x}(k+1)\notag\\
&=\bm{x}(k+1)-\bm{b}^{x}(k+1)\notag\\
&=\bm{x}(k+1)-s(k)Q\left[\frac{1}{s(k)}(\bm{x}(k+1)-\bm{b}^{x}(k))\right]-\bm{b}^{x}(k)\notag\\
&=s(k)(\theta^{x}(k)-Q[\theta^{x}(k)]).
\end{align}

Similarly, by subtracting $\bm{b}^{u}(k)$ from both sides of (\ref{CptAlg2B}), it can be obtained that
\begin{align}\label{YB2}
&\bm{u}(k+1)-\bm{b}^{u}(k)\notag\\
&=\bm{u}(k)-\beta\bm{L}\bm{u}(k)+\beta\bm{L}\bm{e}^{u}(k)+\bm{g}(k+1)-\bm{g}(k)-\bm{b}^{u}(k)\notag\\
&=(\bm{\mathrm{I}}_{nm}+\beta\bm{L})\bm{e}^{u}(k)-\beta\bm{L}\bm{u}(k)+\bm{g}(k+1)-\bm{g}(k)\notag\\
&=s(k)((\bm{\mathrm{I}}_{nm}+\beta\bm{L})\LX{\frac{\bm{e}^{u}(k)}{s(k)}}-\frac{1}{s(k)}(\beta\bm{L}\bm{u}(k)\notag\\
&\quad-\bm{g}(k+1)+\bm{g}(k)))\notag\\
&=s(k)\theta^{u}(k),
\end{align}

where
$$
\theta^{u}(k)=(\bm{\mathrm{I}}_{nm}+\beta\bm{L})\LX{\frac{\bm{e}^{u}(k)}{s(k)}}-\frac{1}{s(k)}(\beta\bm{L}\bm{u}(k)\LX{-\bm{g}(k+1)+\bm{g}(k)}).
$$
Based on (\ref{Estimation2a2}) and (\ref{YB2}), we have
\begin{align}\label{ErrorY2}
&\bm{e}^{u}(k+1)\notag\\
&=\bm{u}(k+1)-\bm{b}^{u}(k+1)\notag\\
&=\bm{u}(k+1)-s(k)Q\left[\frac{1}{s(k)}(\bm{u}(k+1)-\bm{b}^{u}(k))\right]-\bm{b}^{u}(k)\notag\\
&=s(k)(\theta^{u}(k)-Q[\theta^{u}(k)]).
\end{align}

The proof of nonsaturation of the uniform quantizer is equivalent to showing that for any $k\geq0$, $\|\theta^{x}(k)\|_{\infty}\leq\LX{\mathcal{K}}+\frac{1}{2}$ and $\|\theta^{u}(k)\|_{\infty}\leq\LX{\mathcal{K}}+\frac{1}{2}$. The proof is based on
induction. We begin by showing the quantizer is not saturated at $k=0$.

Note that $\frac{\bm{e}^{x}(0)}{s(0)}=\frac{\bm{x}(0)-\bm{b}^{x}(0)}{s(0)}=\frac{\bm{x}(0)}{s(0)}$, $\frac{\bm{e}^{u}(0)}{s(0)}=\frac{\bm{u}(0)-\bm{b}^{u}(0)}{s(0)}=\frac{\bm{u}(0)}{s(0)}$ and $\|\bm{x}(0)\|_{\infty}\leq C_{x}$, $\|\bm{u}(0)\|_{\infty}\leq C_{u}$, we have
\begin{align}\label{theta0x}
&\|\theta^{x}(0)\|_{\infty}\notag\\
&=\|(\bm{\mathrm{I}}_{nm}+\beta\bm{L})\LX{\frac{\bm{e}^{x}(0)}{s(0)}}-\frac{1}{s(0)}(\beta\bm{L}\bm{x}(0)+\delta\bm{u}(0))\|_{\infty}\notag\\
&\leq\|(\bm{\mathrm{I}}_{nm}+\beta\bm{L})\frac{\bm{x}(0)}{s(0)}-\frac{\beta}{s(0)}\bm{L}\bm{x}(0)\|_{\infty}+\|\frac{\delta}{s(0)}\bm{u}(0)\|_{\infty}\notag\\
&\leq\frac{C_{x}+\delta C_{u}}{s(0)}\leq\mathcal{K}+\frac{1}{2},
\end{align}
where the last inequality can be guaranteed by the condition in (\ref{ScaF2}), and
\begin{align}\label{theta0y}
&\|\theta^{u}(0)\|_{\infty}\notag\\
&=\|(\bm{\mathrm{I}}_{nm}+\beta\bm{L})\LX{\frac{\bm{e}^{u}(0)}{s(0)}}-\frac{1}{s(0)}(\beta\bm{L}\bm{u}(0)-\bm{g}(1)\notag\\
&\quad+\bm{g}(0))\|_{\infty}\notag\\
&=\|(\bm{\mathrm{I}}_{nm}+\beta\bm{L})\frac{\bm{u}(0)}{s(0)}-\frac{1}{s(0)}(\beta\bm{L}\bm{u}(0)+\bm{g}(0))\notag\\
&\quad+\frac{1}{s(0)}\bm{g}(1)\|_{\infty}\notag\\
&=\|\frac{1}{s(0)}\nabla F(\bm{x}(1))\|_{\infty}\notag\\
&=\|\frac{1}{s(0)}\nabla F(\bm{x}(0)-\delta\bm{u}(0))\|_{\infty}\leq\mathcal{K}+\frac{1}{2},
\end{align}
where the third equality holds due to $\bm{u}(0)=\nabla F(\bm{x}(0))$;
and the last inequality can be guaranteed by the condition in (\ref{ScaF2}). Hence, the quantizer is unsaturated at $k=0$. Now, we assume that the quantizer is not saturated at $\LX{k=0,\dots,p}$. Then, by (\ref{Error2}) and \LX{(\ref{ErrorY2})}, it can be calculated that
\begin{subequations}\label{EpsilonK2}
\begin{align}
\LX{\sup_{1\leq k\leq p+1}\|\frac{\bm{e}^{x}(k)}{s(k)}\|_{\infty}\leq\frac{1}{2\mu},}\label{EpsilonK2a}\\
\LX{\sup_{1\leq k\leq p+1}\|\frac{\bm{e}^{u}(k)}{s(k)}\|_{\infty}\leq\frac{1}{2\mu}.}\label{EpsilonK2b}
\end{align}
\end{subequations}

We proceed to show that the quantizer is unsaturated for $k=p+1$.
\begin{align}\label{thetaP1a2}
&\|\theta^{x}(p+1)\|_{\infty}\notag\\
&=\|(\bm{\mathrm{I}}_{nm}+\beta\bm{L})\LX{\frac{\bm{e}^{x}(p+1)}{s(p+1)}}-\frac{1}{s(p+1)}(\beta\bm{L}\bm{x}(p+1)\notag\\
&\quad+\delta\bm{u}(p+1))\|_{\infty}\notag\\
&\leq\frac{(1+2\beta d)}{2\mu}+\frac{1}{s(p+1)}\|\beta\bm{L}\bm{x}(p+1)+\delta\bm{u}(p+1)\|_{\infty}\notag\\
&\leq\frac{(1+2\beta d)}{2\mu}+\frac{1}{s(p+1)}\|\beta\bm{L}\bm{x}(p+1)+\delta\bm{u}(p+1)\|\notag\\
&\leq\frac{(1+2\beta d)}{2\mu}+\frac{\sqrt{3}}{s(p+1)}(\|\beta\bm{L}\bm{x}(p+1)\|^{2}\notag\\
&\quad+\delta^{2}\|\bm{u}(p+1)-\bar{\bm{u}}(p+1)\|^{2}+\delta^{2}\|\bar{\bm{u}}(p+1)\|^{2})^{\frac{1}{2}}\notag\\
&\leq\frac{(1+2\beta d)}{2\mu}+\frac{\sqrt{3}}{s(p+1)}(\beta^{2}\rho^{2}(L)\|\bm{x}(p+1)-\bar{\bm{x}}(p+1)\|^{2}\notag\\
&\quad+\delta^{2}\|\bm{u}(p+1)-\bar{\bm{u}}(p+1)\|^{2}+\delta^{2}\|\bar{\bm{g}}(p+1)\|^{2})^{\frac{1}{2}}\notag\\
&\leq\frac{(1+2\beta d)}{2\mu}+\frac{\sigma_{6}}{s(p+1)}(\|\bm{x}(p+1)-\bar{\bm{x}}(p+1)\|^{2}\notag\\
&\quad+\|\bm{u}(p+1)-\bar{\bm{u}}(p+1)\|^{2}+\|\bar{\bm{g}}(p+1)\|^{2})^{\frac{1}{2}},
\end{align}
where the first inequality holds due to \LX{(\ref{EpsilonK2a})} and $\|\bm{\mathrm{I}}_{nm}+\beta\bm{L}\|_{\infty}=1+2\beta d$; the third inequality holds due to the Cauchy--Schwarz inequality; and the last inequality holds due to (\ref{Lemma1Ab}) \LX{and} (\ref{barU2}).

Similarly, we have
\blue{
\begin{align}\label{thetaP1u2}
&\|\theta^{u}(p+1)\|_{\infty}\notag\\
&=\|(\bm{\mathrm{I}}_{nm}+\beta\bm{L})\LX{\frac{\bm{e}^{u}(p+1)}{s(p+1)}}-\frac{1}{s(p+1)}(\beta\bm{L}\bm{u}(p+1)\notag\\
&\quad-\bm{g}(p+2)+\bm{g}(p+1))\|_{\infty}\notag\\
&\leq\frac{(1+2\beta d)}{2\mu}+\frac{1}{s(p+1)}\|\beta\bm{L}\bm{u}(p+1)-\bm{g}(p+2)\notag\\
&\quad+\bm{g}(p+1)\|_{\infty}\notag\\
&\leq\frac{(1+2\beta d)}{2\mu}+\frac{1}{s(p+1)}(\|\beta\bm{L}\bm{u}(p+1)\|\notag\\
&\quad+\|\bm{g}(p+2)-\bm{g}(p+1)\|)\notag\\
&\leq\frac{(1+2\beta d)}{2\mu}+\frac{1}{s(p+1)}(\beta\rho(L)\|\bm{u}(p+1)-\bar{\bm{u}}(p+1)\|\notag\\
&\quad+L_{f}(\beta\rho(L)\|\bm{x}(p+1)-\bar{\bm{x}}(p+1)\|+\beta\rho(L)\|\bm{e}^{x}(p+1)\|\notag\\
&\quad+\delta\|\bm{u}(p+1)-\bar{\bm{u}}(p+1)\|+\delta\|\bar{\bm{g}}(p+1)\|))\notag\\
&\leq\frac{(1+2\beta d)}{2\mu}+\frac{\sqrt{nm}L_{f}\beta\rho(L)}{2\mu}\notag\\
&\quad+\frac{1}{s(p+1)}((\beta\rho(L)+L_{f}\delta)\|\bm{u}(p+1)-\bar{\bm{u}}(p+1)\|\notag\\
&\quad+L_{f}(\beta\rho(L)\|\bm{x}(p+1)-\bar{\bm{x}}(p+1)\|+\delta\|\bar{\bm{g}}(p+1)\|)),
%
%
\end{align}
where the first inequality holds due to \LX{(\ref{EpsilonK2b})};
the third inequality holds due to (\ref{Lemma1Ab}), (\ref{ErrorgReg}) and $\rho(\bm{K})\leq1$; the last inequality holds due to $\|x\|\leq\sqrt{nm}\|x\|_{\infty}$ and (\ref{EpsilonK2a}).}

Then, from (\ref{ErrorX2}), \LX{(\ref{OptimalF2})--(\ref{ErrorU21})}, we have
\begin{equation}\label{Lambdak}
\Lambda(k+1)\preceq\Phi\Lambda(k)+\Gamma(k),
\end{equation}
where $\Lambda(k)=[\|\bm{x}(k)-\bar{\bm{x}}(k)\|^{2},~\|\bm{u}(k)-\bar{\bm{u}}(k)\|^{2},~n(f(\bar{x}(k))-f^{\star})]^{T}$, $\Phi$ is defined in (\ref{PPhi}) and
\begin{align}\label{Gamma}
\Gamma(k)=\sigma_{5}\left[
  \begin{array}{c}
   \|e^{x}(k)\|^{2} \\
    \|e^{u}(k)\|^{2}+4L_{f}^{2}\|e^{x}(k)\|^{2} \\
    0 \\
  \end{array}
\right].
\end{align}

Based on (\ref{EpsilonK2}) and (\ref{Gamma}), it can be obtained that
\LX{\begin{equation}\label{Gammak}
\|\Gamma(k)\|\leq\frac{1}{4\mu^{2}}\sigma_{2}s^{2}(k).
\end{equation}}

Based on (\ref{Lambdak}), one can have
\begin{equation}\label{Lambdak1}
\|\Lambda(k+1)\|\leq\|\Phi^{k+1}\Lambda(0)\|+\sum_{\tau=0}^{k}\|\Phi^{\tau}\|\|\Gamma(k-\tau)\|.
\end{equation}


From \citet[Corollary 8.1.29]{Horn2012Matrix} and Lemma~\ref{lem-1}, we know that $\rho(\Phi)\leq\LX{\bar{\rho}}$. Hence, based on (\ref{hatrho}), it can be derived that
\begin{equation}\label{barrho}
\|\Phi^{k}\|\leq h\bar{\rho}^{k}.
\end{equation}
Then, \LX{for $\forall k=0,\dots,p$,} based on \LX{(\ref{Gammak})--(\ref{barrho})}, we have
\begin{equation}\label{Lambdak2}
\|\Lambda(k+1)\|\leq h\bar{\rho}^{k+1}\|\Lambda(0)\|+\frac{h}{4\mu^{2}}\sigma_{2}\sum_{\tau=0}^{k}\bar{\rho}^{\tau}s^{2}(k-\tau).
\end{equation}

From (\ref{OptimalF2c}), we have
\begin{equation}\label{OptimalF2d}
\|\bar{\bm{g}}(k)\|^{2}\leq\frac{\LX{4}\sigma_{7}}{\delta(1-2 \delta L_{f})}\|\Lambda(k)\|.
\end{equation}

Then, from $\mu\in(\sqrt{\bar{\rho}},1)$, it can be guaranteed that $\mu^{2}>\bar{\rho}$.

Based on (\ref{K-level2}) and (\ref{ScaF2}), we know that $s(0)\geq\sqrt{\frac{4\|\Lambda(0)\|\mu^{2}(\mu^{2}-\bar{\rho})}{\sigma_{2}}}$ and $\mathcal{K}\geq\vartheta_{1}$, it can be guaranteed that
\begin{align}\label{thetaP1a3a}
&\|\theta^{x}(p+1)\|_{\infty}\notag\\
&\leq\frac{(1+2\beta d)}{2\mu}+\frac{\sigma_{6}}{s(p+1)}(\|\bm{x}(p+1)-\bar{\bm{x}}(p+1)\|^{2}\notag\\
&\quad+\|\bm{u}(p+1)-\bar{\bm{u}}(p+1)\|^{2}+\|\bar{\bm{g}}(p+1)\|^{2})^{\frac{1}{2}}\notag\\
&\leq\frac{(1+2\beta d)}{2\mu}+\frac{\sigma_{6}}{s(p+1)}\sqrt{(2+\frac{\LX{4}\sigma_{7}}{\delta(1-2 \delta L_{f})})}\|\Lambda(p+1)\|^{\frac{1}{2}}\notag\\
&=\frac{(1+2\beta d)}{2\mu}+\frac{\sigma_{3}}{s(p+1)}\|\Lambda(p+1)\|^{\frac{1}{2}}\notag\\
&\leq\frac{\sigma_{3}}{s(p+1)}(h\bar{\rho}^{p+1}\|\Lambda(0)\|+\frac{h}{4\mu^{2}}\sigma_{2}\sum_{\tau=0}^{p}\bar{\rho}^{\tau}s^{2}(p-\tau))^{\frac{1}{2}}\notag\\
&\quad+\frac{(1+2\beta d)}{2\mu}\notag\\
&=\frac{\sigma_{3}}{s(p+1)}(h\bar{\rho}^{p+1}\|\Lambda(0)\|+\frac{h}{4\mu^{2}}\sigma_{2}s^{2}(p)\sum_{\tau=0}^{p}(\frac{\bar{\rho}}{\mu^{2}})^{\tau})^{\frac{1}{2}}\notag\\
&\quad+\frac{(1+2\beta d)}{2\mu}\notag\\
&=\sigma_{3}(\frac{h\|\Lambda(0)\|}{s^{2}(0)}(\frac{\bar{\rho}}{\mu^{2}})^{p+1}+\frac{h\sigma_{2}}{4\mu^{2}(\mu^{2}-\bar{\rho})}(1-(\frac{\bar{\rho}}{\mu^{2}})^{p+1}))^{\frac{1}{2}}\notag\\
&\quad+\frac{(1+2\beta d)}{2\mu}\notag\\
&\leq\LX{\vartheta_1+\frac{1}{2}}\leq\mathcal{K}+\frac{1}{2},
\end{align}
where the first inequality holds due to (\ref{thetaP1a2}); the second inequality holds due to (\ref{OptimalF2d}); the third inequality holds due to (\ref{Lambdak2}).

Based on (\ref{K-level2}) and (\ref{ScaF2}), we know that $s(0)\geq\sqrt{\frac{4\|\Lambda(0)\|\mu^{2}(\mu^{2}-\bar{\rho})}{\sigma_{2}}}$ and $\mathcal{K}\geq\vartheta_{2}$, it can be guaranteed that
\begin{align}\label{thetaP1a3b}
&\|\theta^{u}(p+1)\|_{\infty}\notag\\
&\leq\frac{(1+2\beta d)}{2\mu}+\frac{\sqrt{nm}L_{f}\beta\rho(L)}{2\mu}\notag\\
&\quad+\frac{1}{s(p+1)}((\beta\rho(L)+L_{f}\delta)\|\Lambda(p+1)\|^{\frac{1}{2}}\notag\\
&\quad+L_{f}(\beta\rho(L)\|\Lambda(p+1)\|^{\frac{1}{2}}+\delta\sqrt{\frac{\LX{4}\sigma_{7}}{\delta(1-2 \delta L_{f})}}\|\Lambda(p+1)\|^{\frac{1}{2}}))\notag\\
&\leq\frac{(1+2\beta d)}{2\mu}+\frac{\sqrt{nm}L_{f}\beta\rho(L)}{2\mu}\notag\\
&\quad+\sigma_{4}(\frac{h\|\Lambda(0)\|}{s^{2}(0)}(\frac{\bar{\rho}}{\mu^{2}})^{p+1}+\frac{h\sigma_{2}}{4\mu^{2}(\mu^{2}-\bar{\rho})}(1-(\frac{\bar{\rho}}{\mu^{2}})^{p+1})^{\frac{1}{2}})\notag\\
&\LX{\leq\vartheta_2+\frac{1}{2}\leq\mathcal{K}+\frac{1}{2},}
\end{align}
\blue{where the first inequality holds due to (\ref{thetaP1u2}) and (\ref{OptimalF2d}).
As a result, when $k=p+1$, the quantizer is also unsaturated. Therefore, by induction, we conclude that the quantizer is never saturated. This proof is complete.}

\section{Proof of Theorem \ref{thm-3}}\label{Appendix-E}

From (\ref{Lambdak2}) and $s(0)\geq\sqrt{\frac{4\|\Lambda(0)\|\mu^{2}(\mu^{2}-\bar{\rho})}{\sigma_{2}}}$, we have
\begin{align}\label{LambdakThm1}
&\|\Lambda(k+1)\|\notag\\
&\leq h\bar{\rho}^{k+1}\|\Lambda(0)\|+\frac{h}{4\mu^{2}}\sigma_{2}\sum_{\tau=0}^{k}\bar{\rho}^{\tau}s^{2}(k-\tau)\notag\\
&= h\bar{\rho}^{k+1}\|\Lambda(0)\|+\frac{h\sigma_{2}s^{2}(k)}{4(\mu^{2}-\bar{\rho})}(1-(\frac{\bar{\rho}}{\mu^{2}})^{p+1})\notag\\
&=\bar{\rho}^{k+1}(h\|\Lambda(0)\|-\frac{h\sigma_{2}s^{2}(0)}{4\mu^{2}(\mu^{2}-\bar{\rho})})+\frac{h\sigma_{2}s^{2}(0)}{4\mu^{2}(\mu^{2}-\bar{\rho})}\mu^{2(k+1)}\notag\\
&\leq\frac{h\sigma_{2}s^{2}(0)}{4\mu^{2}(\mu^{2}-\bar{\rho})}\mu^{2(k+1)}.
\end{align}
This proof is complete.
%
%
%

\section{Proof of Theorem \ref{thm-4}}\label{Appendix-F}

\blue{Note that
\begin{equation}\label{eta02a}
\lim_{\beta\rightarrow0}\sigma_{3}\sqrt{\frac{h\sigma_{2}}{4\mu^{2}(\mu^{2}-\bar{\rho})}}+\frac{(1+2\beta d)}{2\mu}=\frac{1}{2\mu},
\end{equation}
and
\begin{equation}\label{eta02b}
\lim_{\beta\rightarrow0}\frac{(1+2\beta d)}{2\mu}+\frac{\sqrt{nm}L_{f}\beta\rho(L)}{2\mu}+\sigma_{4}\sqrt{\frac{h\sigma_{2}}{4\mu^{2}(\mu^{2}-\bar{\rho})}}=\frac{1}{2\mu}.
\end{equation}

Then, for any given $\mathcal{K}\geq1$, there exists
$ \beta^{\star}\in(0,\frac{\sqrt{2}}{2\rho(L)})$
and $\delta^{\star}\in(0,\min\{~\frac{\sqrt{c_{1}\Theta_{1}}}{2},~\frac{1}{4L_{f}},~\frac{2}{\nu},~\frac{1}{8+2L_{f}},~\frac{\sqrt{c_{1}}}{8L_{f}}\})$
such that
\begin{equation}\label{eta02c}
\lim_{\mu\rightarrow1}\sigma_{3}\sqrt{\frac{h\sigma_{2}}{4\mu^{2}(\mu^{2}-\bar{\rho})}}+\frac{(1+2\beta^{\star}d)}{2\mu}\leq\mathcal{K}+\frac{1}{2},
\end{equation}
and
\begin{align}\label{eta02d}
\lim_{\mu\rightarrow1}\frac{(1+2\beta d)}{2\mu}&+\frac{\sqrt{nm}L_{f}\beta\rho(L)}{2\mu}\notag\\
&+\sigma_{4}\sqrt{\frac{h\sigma_{2}}{4\mu^{2}(\mu^{2}-\bar{\rho})}}\leq\mathcal{K}+\frac{1}{2}.
\end{align}}

Hence, there exists $\mu^{\star}\in(\sqrt{\bar{\rho}},1)$ such that $\vartheta_1\leq\mathcal{K}$ and $\vartheta_2\leq\mathcal{K}$. Thus, $(\mu^{\star},\beta^{\star},\delta^{\star})\in\Pi$, and hence $\Pi$ is nonempty, where $\Pi$ is defined in Theorem~\ref{thm-4}. The proof of the convergence result is similar to that of Theorem~\ref{thm-3}.

\section{Proof of Proposition \ref{prop}} \label{Appendix-A}
\LX{For simplicity of the notation, define $b_{j}(k)=b_{j}^{x}(k)$, $e_{j}(k)=e_{j}^{x}(k)$, $\bm{e}(k)=[e^T_{1}(k),\dots,e^T_{n}(k)]^{T}$.
Then, the update rule (\ref{Alg}) can be rewritten in a compact form:}
  \begin{subequations}\label{CptAlg2}
  \begin{align}
  \bm{x}(k+1)&=\bm{x}(k)-\xi\bm{L}\bm{x}(k)-\varphi\bm{u}(k)-\sigma\nabla F(\bm{x}(k))\notag\\
  &\quad+\xi\bm{L}\bm{e}(k),\label{CptAlgA}\\
  \quad\quad~\bm{u}(k+1)&=\bm{u}(k)+\varphi\bm{L}\bm{x}(k)-\varphi\bm{L}\bm{e}(k),\notag\\
  &\quad~\forall\bm{x}(0)\in\mathbb{R}^{nm},~\sum_{j=1}^{n}u_{j}(0)=\bm{0}_{m}.\label{CptAlgB}
  \end{align}
  \end{subequations}


Noting that $\nabla F$ is Lipschitz-continuous with constant $L_{f}>0$ as assumed in Assumption \ref{assum-2}, we have
\begin{align}\label{LipIne}
&\|\bm{g}(k+1)-\bm{g}(k)\|^{2}\notag\\
&\leq L_{f}^{2}\|\bm{x}(k+1)-\bm{x}(k)\|^{2}\notag\\
&=L_{f}^{2}\|-\xi\bm{L}\bm{x}(k)-\varphi\bm{u}(k)-\sigma\nabla F(\bm{x}(k))+\xi\bm{L}\bm{e}(k)\|^{2}\notag\\
&\leq3L_{f}^{2}(\xi^{2}\bm{x}^{T}(k)\bm{L}^{2}\bm{x}(k)+\varphi^{2}\|\bm{u}(k)+\frac{\sigma}{\varphi}\bm{g}(k)\|^{2}\notag\\
&\quad+\xi^{2}\bm{e}^{T}(k)\bm{L}^{2}\bm{e}(k)),
\end{align}
where the second inequality holds due to the Cauchy--Schwarz inequality.
Note that $\bar{\bm{u}}(k)=\bm{1}_{n}\otimes\bar{u}(k)$ and $\bar{u}(k)=\frac{1}{n}(\bm{1}_{n}^{T}\otimes\bm{\mathrm{I}}_{m})\bm{u}(k)$, then from (\ref{CptAlgB}), we know that $\bar{u}(k+1)=\bar{u}(k)$, and due to $\sum_{j=1}^{n}u_{j}(0)=\bm{0}_{m}$, we have
\begin{equation}\label{barU}
\bar{u}(k)\equiv\bm{0}_{m}.
\end{equation}

From (\ref{CptAlgA}), it can be calculated that
  \begin{equation}\label{barX}
  \bar{\bm{x}}(k+1)=\bar{\bm{x}}(k)-\sigma\bar{\bm{g}}(k).
  \end{equation}

Denote $\bm{b}(k)=[b^T_{1}(k),\dots,b^T_{n}(k)]^{T}$. Based on (\ref{Encoder2Ba}), we know that
\begin{align}\label{Estimation}
\bm{b}(k)&=s(k-1)Q\left[\frac{1}{s(k-1)}(\bm{x}(k)-\bm{b}(k-1))\right]\notag\\
&\quad+\bm{b}(k-1).
\end{align}

By subtracting $\bm{b}(k)$ from both sides of (\ref{CptAlgA}), it can be obtained that
\begin{align}\label{XB}
&\bm{x}(k+1)-\bm{b}(k)\notag\\
&=\bm{x}(k)-\xi\bm{L}\bm{x}(k)-\varphi\bm{u}(k)-\sigma\bm{g}(k)+\xi\bm{L}\bm{e}(k)-\bm{b}(k)\notag\\
&=(\bm{\mathrm{I}}_{nm}+\xi\bm{L})\bm{e}(k)-\xi\bm{L}\bm{x}(k)-\varphi(\bm{u}(k)+\frac{\sigma}{\varphi}\bm{g}(k))\notag\\
&=s(k)\left[(\bm{\mathrm{I}}_{nm}+\xi\bm{L})\LX{\frac{\bm{e}(k)}{s(k)}}-\frac{1}{s(k)}(\xi\bm{L}\bm{x}(k)\right.\notag\\
&\quad\left.+\varphi(\bm{u}(k)+\frac{\sigma}{\varphi}\bm{g}(k)))\right]\notag\\
&=s(k)\theta(k),
\end{align}
where
$$
\theta(k)=(\bm{\mathrm{I}}_{nm}+\xi\bm{L})\LX{\frac{\bm{e}(k)}{s(k)}}-\frac{1}{s(k)}(\xi\bm{L}\bm{x}(k)+\varphi(\bm{u}(k)+\frac{\sigma}{\varphi}\bm{g}(k))).
$$

Based on (\ref{Estimation}) and (\ref{XB}), we have
\begin{align}\label{Error}
&\bm{e}(k+1)\notag\\
&=\bm{x}(k+1)-\bm{b}(k+1)\notag\\
&=\bm{x}(k+1)-s(k)Q\left[\frac{1}{s(k)}(\bm{x}(k+1)-\bm{b}(k))\right]-\bm{b}(k)\notag\\
&=s(k)(\theta(k)-Q[\theta(k)]).
\end{align}

The proof of nonsaturation of the uniform quantizer is equivalent to showing that for any $k\geq0$, $\|\theta(k)\|_{\infty}\leq \LX{\mathcal{K}}+\frac{1}{2}$. The proof is based on
induction. We begin by showing the quantizer is not saturated at $k=0$.

Note that $\frac{\bm{e}(0)}{s(0)}=\frac{\bm{x}(0)-\bm{b}(0)}{s(0)}=\frac{\bm{x}(0)}{s(0)}$, and $\|\bm{x}(0)\|_{\infty}\leq C_{x}$, $\|\bm{u}(0)\|_{\infty}\leq C_{u}$, $\|\bm{g}(0)\|_{\infty}\leq C_{g}$, we have
\begin{align}\label{theta0}
&\|\theta(0)\|_{\infty}\notag\\
&=\|(\bm{\mathrm{I}}_{nm}+\xi\bm{L})\LX{\frac{\bm{e}(0)}{s(0)}}-\frac{1}{s(0)}(\xi\bm{L}\bm{x}(0)\notag\\
&\quad+\varphi(\bm{u}(0)+\frac{\sigma}{\varphi}\bm{g}(0)))\|_{\infty}\notag\\
&\leq\|(\bm{\mathrm{I}}_{nm}+\xi\bm{L})\frac{\bm{x}(0)}{s(0)}-\frac{\xi}{s(0)}\bm{L}\bm{x}(0)\|_{\infty}\notag\\
&\quad+\|\frac{\varphi}{s(0)}\bm{u}(0)\|_{\infty}+\|\frac{\sigma}{s(0)}\bm{g}(0)\|_{\infty}\notag\\
&\leq\frac{C_{x}+\varphi C_{u}+\sigma C_{g}}{s(0)}\leq\mathcal{K}+\frac{1}{2},
\end{align}
where the last inequality can be guaranteed by the condition in (\ref{ScaF}). Hence, the quantizer is unsaturated at $k=0$. Now, we assume that the quantizer is not saturated at $k=0,\dots,p$. Then, by (\ref{Error}), it can be calculated that
\begin{equation}\label{EpsilonK}
\sup_{1\leq k\leq p+1}\|\frac{\bm{e}(k)}{s(k)}\|_{\infty}\leq\frac{1}{2\mu}.
\end{equation}

We proceed to show that the quantizer is unsaturated for $k=p+1$.
\begin{align}\label{thetaP1a}
&\|\theta(p+1)\|_{\infty}\notag\\
&=\|(\bm{\mathrm{I}}_{nm}+\xi\bm{L})\LX{\frac{\bm{e}(p+1)}{s(p+1)}}-\frac{1}{s(p+1)}(\xi\bm{L}\bm{x}(p+1)\notag\\
&\quad+\varphi(\bm{u}(p+1)+\frac{\sigma}{\varphi}\bm{g}(p+1)))\|_{\infty}\notag\\
&\leq\|(\bm{\mathrm{I}}_{nm}+\xi\bm{L})\LX{\frac{\bm{e}(p+1)}{s(p+1)}}\|_{\infty}+\|\frac{1}{s(p+1)}(\xi\bm{L}\bm{x}(p+1)\notag\\
&\quad+\varphi(\bm{u}(p+1)+\frac{\sigma}{\varphi}\bm{g}(p+1)))\|_{\infty}\notag\\
&\leq\frac{(1+2\xi d)}{2\mu}+\|\frac{1}{s(p+1)}(\xi\bm{L}\bm{x}(p+1)+\varphi(\bm{u}(p+1)\notag\\
&\quad+\frac{\sigma}{\varphi}\bm{g}(p+1)))\|_{\infty}\notag\\
&\leq\frac{(1+2\xi d)}{2\mu}+\|\frac{1}{s(p+1)}(\xi\bm{L}\bm{x}(p+1)+\varphi(\bm{u}(p+1)\notag\\
&\quad+\frac{\sigma}{\varphi}\bm{g}(p+1)))\|\notag\\
&=\frac{(1+2\xi d)}{2\mu}+\frac{1}{s(p+1)}(\|\xi\bm{L}\bm{x}(p+1)\|^{2}\notag\\
&\quad+\|\varphi(\bm{u}(p+1)+\frac{\sigma}{\varphi}\bm{g}(p+1))\|^{2}\notag\\
&\quad+2\xi\varphi\bm{x}^{T}(p+1)\bm{L}(\bm{u}(p+1)+\frac{\sigma}{\varphi}\bm{g}(p+1)))^{\frac{1}{2}}\notag\\
&\leq\frac{(1+2\xi d)}{2\mu}+\frac{1}{s(p+1)}(\xi^{2}\rho^{2}(L)\bm{x}^{T}(p+1)\bm{K}\bm{x}(p+1)\notag\\
&\quad+\frac{\varphi^{3}\rho(L)}{\varphi+\xi}(\bm{u}(p+1)+\frac{\sigma}{\varphi}\bm{g}(p+1))^{T}(\frac{\varphi+\xi}{\varphi}\bm{P})\notag\\
&\quad\times(\bm{u}(p+1)+\frac{\sigma}{\varphi}\bm{g}(p+1))\notag\\
&\quad+2\xi\varphi\rho^{2}(L)\bm{x}^{T}(p+1)\bm{K}\bm{P}(\bm{u}(p+1)+\frac{\sigma}{\varphi}\bm{g}(p+1)))^{\frac{1}{2}}\notag\\
&\leq\frac{(1+2\xi d)}{2\mu}+\frac{1}{s(p+1)}\epsilon_{1}(\bm{x}^{T}(p+1)\bm{K}\bm{x}(p+1)\notag\\
&\quad+(\bm{u}(p+1)+\frac{\sigma}{\varphi}\bm{g}(p+1))^{T}(\frac{\varphi+\xi}{\varphi}\bm{P})\notag\\
&\quad\times(\bm{u}(p+1)+\frac{\sigma}{\varphi}\bm{g}(p+1))\notag\\
&\quad+2\bm{x}^{T}(p+1)\bm{K}\bm{P}(\bm{u}(p+1)+\frac{\sigma}{\varphi}\bm{g}(p+1)))^{\frac{1}{2}},
\end{align}
where the second inequality holds due to (\ref{EpsilonK}); the fourth inequality holds due to (\ref{Lemma1Bb}) and (\ref{Lemma1Bc}).


For the first term in $V(k)$, we know that
\begin{align}\label{Firstterm}
&\bm{x}^{T}(p+1)\bm{K}\bm{x}(p+1)\notag\\
&=\bm{x}^{T}(p)\bm{K}\bm{x}(p)-2\xi\bm{x}^{T}(p)\bm{L}(\bm{x}(p)-\bm{e}(p))\notag\\
&\quad-2\varphi(\bm{x}^{T}(p)-\xi(\bm{x}(p)-\bm{e}(p))^{T}\bm{L})\bm{K}(\bm{u}(p)+\frac{\sigma}{\varphi}\bm{g}(p))\notag\\
&\quad+\varphi^{2}(\bm{u}(p)+\frac{\sigma}{\varphi}\bm{g}(p))^{T}\bm{K}(\bm{u}(p)+\frac{\sigma}{\varphi}\bm{g}(p))\notag\\
&\quad+\xi^{2}(\bm{x}(p)-\bm{e}(p))^{T}\bm{L}^{2}(\bm{x}(p)-\bm{e}(p))\notag\\
&\leq\bm{x}^{T}(p)\bm{K}\bm{x}(p)-2\xi\bm{x}^{T}(p)\bm{L}\bm{x}(p)+\xi\bm{x}^{T}(p)\bm{L}\bm{x}(p)\notag\\
&\quad+\xi\bm{e}^{T}(p)\bm{L}\bm{e}(p)+\xi^{2}(\bm{x}(p)-\bm{e}(p))^{T}\bm{L}^{2}(\bm{x}(p)-\bm{e}(p))\notag\\
&\quad-2\varphi\bm{x}^{T}(p)\bm{K}(\bm{u}(p)+\frac{\sigma}{\varphi}\bm{g}(p))\notag\\
&\quad+\xi^{2}(\bm{x}(p)-\bm{e}(p))^{T}\bm{L}^{2}(\bm{x}(p)-\bm{e}(p))\notag\\
&\quad+\varphi^{2}\|\bm{u}(p)+\frac{\sigma}{\varphi}\bm{g}(p)\|^{2}\notag\\
&\quad+\varphi^{2}(\bm{u}(p)+\frac{\sigma}{\varphi}\bm{g}(p))^{T}\bm{K}(\bm{u}(p)+\frac{\sigma}{\varphi}\bm{g}(p))\notag\\
&\leq\bm{x}^{T}(p)\bm{K}\bm{x}(p)-\xi\bm{x}^{T}(p)\bm{L}\bm{x}(p)+\xi\bm{e}^{T}(p)\bm{L}\bm{e}(p)\notag\\
&\quad+4\xi^{2}\bm{x}^{T}(p)\bm{L}^{2}\bm{x}(p)+4\xi^{2}\bm{e}^{T}(p)\bm{L}^{2}\bm{e}(p)\notag\\
&\quad-2\varphi\bm{x}^{T}(p)\bm{K}(\bm{u}(p)+\frac{\sigma}{\varphi}\bm{g}(p))\notag\\
&\quad+2\varphi^{2}\|\bm{u}(p)+\frac{\sigma}{\varphi}\bm{g}(p)\|^{2}\notag\\
&=\bm{x}^{T}(p)\bm{K}\bm{x}(p)-\bm{x}^{T}(p)(\xi\bm{L}-4\xi^{2}\bm{L}^{2})\bm{x}(p)\notag\\
&\quad+\bm{e}^{T}(p)(\xi\bm{L}+4\xi^{2}\bm{L}^{2})\bm{e}(p)\notag\\
&\quad-2\varphi(\bm{x}(p)-\bm{e}(p)+\bm{e}(p))^{T}\bm{K}(\bm{u}(p)+\frac{\sigma}{\varphi}\bm{g}(p))\notag\\
&\quad+2\varphi^{2}\|\bm{u}(p)+\frac{\sigma}{\varphi}\bm{g}(p)\|^{2}\notag\\
&\leq\bm{x}^{T}(p)\bm{K}\bm{x}(p)-\bm{x}^{T}(p)(\xi\bm{L}-4\xi^{2}\bm{L}^{2})\bm{x}(p)\notag\\
&\quad+\bm{e}^{T}(p)(\xi\bm{L}+4\xi^{2}\bm{L}^{2}+2\varphi\rho(L)\bm{K})\bm{e}(p)\notag\\
&\quad-2\varphi(\bm{x}(p)-\bm{e}(p))^{T}\bm{K}(\bm{u}(p)+\frac{\sigma}{\varphi}\bm{g}(p))\notag\\
&\quad+(\bm{u}(p)+\frac{\sigma}{\varphi}\bm{g}(p))^{T}((\frac{\varphi}{2}+2\varphi^{2}\rho(L))\bm{P})\notag\\
&\quad\times(\bm{u}(p)+\frac{\sigma}{\varphi}\bm{g}(p)),
\end{align}
where the first equality holds due to (\ref{CptAlgA}); the first inequality holds due to the Cauchy--Schwarz inequality; the second inequality holds due to $\rho(\bm{K})=1$; the last inequality holds due to the Cauchy--Schwarz inequality.

For the second term in $V(k)$, it can be calculated that
\begin{align}\label{Secondterm}
&(\bm{u}(p+1)+\frac{\sigma}{\varphi}\bm{g}(p+1))^{T}(\frac{\varphi+\xi}{\varphi}\bm{P})(\bm{u}(p+1)+\frac{\sigma}{\varphi}\bm{g}(p+1))\notag\\
&=(\bm{u}(p)+\frac{\sigma}{\varphi}\bm{g}(p))^{T}(\frac{\varphi+\xi}{\varphi}\bm{P})(\bm{u}(p)+\frac{\sigma}{\varphi}\bm{g}(p))\notag\\
&\quad+2(\xi+\varphi)(\bm{x}(p)-\bm{e}(p))^{T}\bm{K}(\bm{u}(p)+\frac{\sigma}{\varphi}\bm{g}(p))\notag\\
&\quad+(\bm{x}(p)-\bm{e}(p))^{T}(\varphi(\xi+\varphi)\bm{L})(\bm{x}(p)-\bm{e}(p))\notag\\
&\quad+\frac{\sigma^{2}}{\varphi^{2}}(\bm{g}(p+1)-\bm{g}(p))^{T}(\frac{\varphi+\xi}{\varphi}\bm{P})(\bm{g}(p+1)-\bm{g}(p))\notag\\
&\quad+\frac{2\sigma}{\varphi}(\bm{u}(p)+\frac{\sigma}{\varphi}\bm{g}(p))^{T}(\bm{P}+\frac{\xi}{\varphi}\bm{P})(\bm{g}(p+1)-\bm{g}(p))\notag\\
&\quad+2\sigma(\bm{x}(p)-\bm{e}(p))^{T}(\bm{K}+\frac{\xi}{\varphi}\bm{K})(\bm{g}(p+1)-\bm{g}(p))\notag\\
&\leq(\bm{u}(p)+\frac{\sigma}{\varphi}\bm{g}(p))^{T}(\frac{\varphi+\xi}{\varphi}\bm{P})(\bm{u}(p)+\frac{\sigma}{\varphi}\bm{g}(p))\notag\\
&\quad+2(\xi+\varphi)(\bm{x}(p)-\bm{e}(p))^{T}\bm{K}(\bm{u}(p)+\frac{\sigma}{\varphi}\bm{g}(p))\notag\\
&\quad+(\bm{x}(p)-\bm{e}(p))^{T}(\varphi(\xi+\varphi)\bm{L})(\bm{x}(p)-\bm{e}(p))\notag\\
&\quad+\frac{\varphi}{2}(\bm{u}(p)+\frac{\sigma}{\varphi}\bm{g}(p))^{T}\bm{P}(\bm{u}(p)+\frac{\sigma}{\varphi}\bm{g}(p))\notag\\
&\quad+(\bm{g}(p+1)-\bm{g}(p))^{T}(\frac{\sigma^{2}(\xi+\varphi)}{\varphi^{3}}\bm{P})(\bm{g}(p+1)-\bm{g}(p))\notag\\
&\quad+(\bm{g}(p+1)-\bm{g}(p))^{T}(\frac{2\sigma^{2}(\xi+\varphi)^{2}}{\varphi^{5}}\bm{P})(\bm{g}(p+1)-\bm{g}(p))\notag\\
&\quad+\sigma^{2}(\bm{x}(p)-\bm{e}(p))^{T}\bm{K}(\bm{x}(p)-\bm{e}(p))\notag\\
&\quad+(\bm{g}(p+1)-\bm{g}(p))^{T}\bm{K}(\bm{g}(p+1)-\bm{g}(p))\notag\\
&\quad+\frac{2\sigma\xi}{\varphi}(\bm{x}(p)-\bm{e}(p))^{T}\bm{K}(\bm{g}(p+1)-\bm{g}(p))\notag\\
&=(\bm{u}(p)+\frac{\sigma}{\varphi}\bm{g}(p))^{T}(\frac{\varphi+\xi}{\varphi}\bm{P})(\bm{u}(p)+\frac{\sigma}{\varphi}\bm{g}(p))\notag\\
&\quad+2(\xi+\varphi)(\bm{x}(p)-\bm{e}(p))^{T}\bm{K}(\bm{u}(p)+\frac{\sigma}{\varphi}\bm{g}(p))\notag\\
&\quad+(\bm{x}(p)-\bm{e}(p))^{T}(\varphi(\xi+\varphi)\bm{L}+\sigma^{2}\bm{K})(\bm{x}(p)-\bm{e}(p))\notag\\
&\quad+(\bm{g}(p+1)-\bm{g}(p))^{T}(\frac{\sigma^{2}(\xi+\varphi)}{\varphi^{3}}\bm{P}+\frac{2\sigma^{2}(\xi+\varphi)^{2}}{\varphi^{5}}\bm{P}\notag\\
&\quad+\bm{K})(\bm{g}(p+1)-\bm{g}(p))\notag\\
&\quad+\frac{\varphi}{2}(\bm{u}(p)+\frac{\sigma}{\varphi}\bm{g}(p))^{T}\bm{P}(\bm{u}(p)+\frac{\sigma}{\varphi}\bm{g}(p))\notag\\
&\quad+\frac{2\sigma\xi}{\varphi}(\bm{x}(p)-\bm{e}(p))^{T}\bm{K}(\bm{g}(p+1)-\bm{g}(p)),
\end{align}
where the first equality holds due to (\ref{CptAlgB}); the inequality holds due to the Cauchy--Schwarz inequality.

For the third term in $V(k)$, it can be obtained that
\begin{align}\label{Thirdterm}
&\bm{x}^{T}(p+1)\bm{K}\bm{P}(\bm{u}(p+1)+\frac{\sigma}{\varphi}\bm{g}(p+1))\notag\\
&=(\bm{x}(p)-\xi\bm{L}\bm{x}(p)-\varphi\bm{u}(p)-\sigma\bm{g}(p)+\xi\bm{L}\bm{e}(p))^{T}\bm{K}\bm{P}\notag\\
&\quad\times(\bm{u}(p)+\frac{\sigma}{\varphi}\bm{g}(p)+\frac{\sigma}{\varphi}\bm{g}(p+1)-\frac{\sigma}{\varphi}\bm{g}(p)\notag\\
&\quad+\varphi\bm{L}\bm{x}(p)-\varphi\bm{L}\bm{e}(p))\notag\\
&=(\bm{x}^{T}(p)\bm{K}\bm{P}-(\xi+\varphi^{2})(\bm{x}(p)-\bm{e}(p))^{T}\bm{K})\notag\\
&\quad\times(\bm{u}(p)+\frac{\sigma}{\varphi}\bm{g}(p))+\varphi\bm{x}^{T}(p)\bm{K}(\bm{x}(p)-\bm{e}(p))\notag\\
&\quad-(\bm{x}(p)-\bm{e}(p))^{T}(\xi\varphi\bm{L})(\bm{x}(p)-\bm{e}(p))\notag\\
&\quad+\frac{\sigma}{\varphi}(\bm{x}^{T}(p)\bm{K}\bm{P}-\xi(\bm{x}(p)-\bm{e}(p)))^{T}\bm{K})(\bm{g}(p+1)-\bm{g}(p))\notag\\
&\quad-(\varphi\bm{u}(p)+\sigma\bm{g}(p)-\sigma\bar{\bm{g}}(p))^{T}\bm{P}(\bm{u}(p)+\frac{\sigma}{\varphi}\bm{g}(p))\notag\\
&\quad-\sigma(\bm{u}(p)+\frac{\sigma}{\varphi}\bm{g}(p))^{T}\bm{K}\bm{P}(\bm{g}(p+1)-\bm{g}(p))\notag\\
&\leq(\bm{x}^{T}(p)\bm{K}\bm{P}-\xi(\bm{x}(p)-\bm{e}(p))^{T}\bm{K})(\bm{u}(p)+\frac{\sigma}{\varphi}\bm{g}(p))\notag\\
&\quad+\frac{\varphi^{2}}{2}(\bm{x}(p)-\bm{e}(p))^{T}\bm{K}(\bm{x}(p)-\bm{e}(p))+\frac{\varphi^{2}}{2}\|\bm{u}(p)\notag\\
&\quad+\frac{\sigma}{\varphi}\bm{g}(p)\|^{2}+\frac{\varphi}{4}\bm{x}^{T}(p)\bm{K}\bm{x}(p)+\frac{\sigma}{2}\bm{x}^{T}(p)\bm{K}\bm{x}(p)\notag\\
&\quad+(\bm{x}(p)-\bm{e}(p))^{T}(\varphi(\bm{K}-\xi\bm{L}))(\bm{x}(p)-\bm{e}(p))\notag\\
&\quad+\frac{\sigma}{2\varphi^{2}}(\bm{g}(p+1)-\bm{g}(p))^{T}\bm{P}^{2}(\bm{g}(p+1)-\bm{g}(p))\notag\\
&\quad-\varphi(\bm{u}(p)+\frac{\sigma}{\varphi}\bm{g}(p))^{T}\bm{P}(\bm{u}(p)+\frac{\sigma}{\varphi}\bm{g}(p))\notag\\
&\quad+\sigma\bar{\bm{g}}^{T}(p)\bm{P}(\bm{u}(p)+\frac{\sigma}{\varphi}\bm{g}(p))\notag\\
&\quad-\frac{\sigma\xi}{\varphi}(\bm{x}(p)-\bm{e}(p))^{T}\bm{K}(\bm{g}(p+1)-\bm{g}(p))\notag\\
&\quad+\sigma^{2}(\bm{u}(p)+\frac{\sigma}{\varphi}\bm{g}(p))^{T}\bm{P}^{2}(\bm{u}(p)+\frac{\sigma}{\varphi}\bm{g}(p))\notag\\
&\quad+\frac{1}{4}\|\bm{g}(p+1)-\bm{g}(p)\|^{2}\notag\\
&\leq\bm{x}^{T}(p)\bm{K}\bm{P}(\bm{u}(p)+\frac{\sigma}{\varphi}\bm{g}(p))\notag\\
&\quad-\xi(\bm{x}(p)-\bm{e}(p))^{T}\bm{K}(\bm{u}(p)+\frac{\sigma}{\varphi}\bm{g}(p))\notag\\
&\quad+(\bm{x}(p)-\bm{e}(p))^{T}(\frac{\varphi^{2}}{2}\bm{K}+\varphi(\bm{K}-\xi\bm{L}))(\bm{x}(p)-\bm{e}(p))\notag\\
&\quad-(\bm{u}(p)+\frac{\sigma}{\varphi}\bm{g}(p))^{T}(\varphi\bm{P}-\frac{\varphi^{2}}{2}\bm{\mathrm{I}}_{nm}-\sigma^{2}\bm{P}^{2}\notag\\
&\quad-4\sigma\bm{P}^{2})(\bm{u}(p)+\frac{\sigma}{\varphi}\bm{g}(p))+\bm{x}^{T}(p)(\frac{\varphi}{4}\bm{K}+\frac{\sigma}{2}\bm{K})\bm{x}(p)\notag\\
&\quad-\frac{\sigma\xi}{\varphi}(\bm{x}(p)-\bm{e}(p))^{T}\bm{K}(\bm{g}(p+1)-\bm{g}(p))\notag\\
&\quad+(\bm{g}(p+1)-\bm{g}(p))^{T}(\frac{\sigma}{2\varphi^{2}}\bm{P}^{2}+\frac{1}{4}\bm{\mathrm{I}}_{nm})\notag\\
&\quad\times(\bm{g}(p+1)-\bm{g}(p))+\frac{\sigma}{16}\|\bar{\bm{g}}(p)\|^{2},
\end{align}
where the first equality holds due to (\ref{CptAlg2}); the second equality holds due to (\ref{Lemma1Aa}) and (\ref{Lemma1Bb});
the first inequality holds due to the Cauchy--Schwarz inequality and $\rho(\bm{K})=1$; the last inequality holds due to the Cauchy--Schwarz inequality.

Then, based on \LX{(\ref{Firstterm})--(\ref{Thirdterm})}, we have
\begin{align}\label{Vk1l}
&V(p+1)\notag\\
&\leq V(p)-\bm{x}^{T}(p)(\xi\bm{L}-4\xi^{2}\bm{L}^{2}-\frac{\varphi}{2}\bm{K}-\sigma\bm{K})\bm{x}(p)\notag\\
&\quad+\bm{e}^{T}(p)(\xi\bm{L}+4\xi^{2}\bm{L}^{2}+2\varphi\rho(L)\bm{K})\bm{e}(p)\notag\\
&\quad-(\bm{u}(p)+\frac{\sigma}{\varphi}\bm{g}(p))^{T}(\varphi\bm{P}-\varphi^{2}\bm{\mathrm{I}}_{nm}-2\sigma^{2}\bm{P}^{2}\notag\\
&\quad-2\varphi^{2}\rho(L)\bm{P}-8\sigma\bm{P}^{2})(\bm{u}(p)+\frac{\sigma}{\varphi}\bm{g}(p))\notag\\
&\quad+(\bm{x}(p)-\bm{e}(p))^{T}(\varphi(\xi+\varphi)\bm{L}+\sigma^{2}\bm{K}+\varphi^{2}\bm{K}\notag\\
&\quad+2\varphi(\bm{K}-\xi\bm{L}))(\bm{x}(p)-\bm{e}(p))\notag\\
&\quad+(\bm{g}(p+1)-\bm{g}(p))^{T}(\frac{\sigma^{2}(\xi+\varphi)}{\varphi^{3}}\bm{P}+\frac{2\sigma^{2}(\xi+\varphi)^{2}}{\varphi^{5}}\bm{P}\notag\\
&\quad+\bm{K}+\frac{\sigma}{\varphi^{2}}\bm{P}^{2}+\frac{1}{2}\bm{\mathrm{I}}_{nm})(\bm{g}(p+1)-\bm{g}(p))+\frac{\sigma}{8}\|\bar{\bm{g}}(p)\|^{2}\notag\\
&\leq V(p)-\bm{x}^{T}(p)(\xi\underline{\rho}(L)\bm{K}-4\xi^{2}\rho^{2}(L)\bm{K}-\frac{\varphi}{2}\bm{K}\notag\\
&\quad-\sigma\bm{K})\bm{x}(p)+(\xi\rho(L)+2\varphi\rho(L)+4\xi^{2}\rho^{2}(L))\|\bm{e}(p)\|^{2}\notag\\
&\quad-(\bm{u}(p)+\frac{\sigma}{\varphi}\bm{g}(p))^{T}(\varphi\bm{P}-\varphi^{2}\rho(L)\bm{P}-\frac{2\sigma^{2}}{\underline{\rho}(L)}\bm{P}\notag\\
&\quad-2\varphi^{2}\rho(L)\bm{P}-\frac{8\sigma}{\underline{\rho}(L)}\bm{P})(\bm{u}(p)+\frac{\sigma}{\varphi}\bm{g}(p))\notag\\
&\quad+2(\varphi(\xi+\varphi)\rho(L)+\sigma^{2}+\varphi^{2}+2\varphi)\notag\\
&\quad\times(\bm{x}^{T}(p)\bm{K}\bm{x}(p)+\|\bm{e}(p)\|^{2})\notag\\
&\quad+3\xi^{2}L_{f}^{2}(\frac{\sigma^{2}(\xi+\varphi)}{\varphi^{3}}\rho(L)+\frac{2\sigma^{2}(\xi+\varphi)^{2}}{\varphi^{5}}\rho(L)\notag\\
&\quad+\rho^{2}(L)+\frac{\sigma}{\varphi^{2}}+\frac{1}{2}\rho^{2}(L))[\bm{x}^{T}(p)\bm{K}\bm{x}(p)+\|\bm{e}(p)\|^{2}]\notag\\
&\quad+3\varphi^{2}L_{f}^{2}(\frac{\sigma^{2}(\xi+\varphi)}{\varphi^{3}}+\frac{2\sigma^{2}(\xi+\varphi)^{2}}{\varphi^{5}}+\frac{\sigma}{\varphi^{2}\underline{\rho}(L)}\notag\\
&\quad+\frac{1}{2}\rho(L)+\rho(L))(\bm{u}(p)+\frac{\sigma}{\varphi}\bm{g}(p))^{T}\bm{P}(\bm{u}(p)+\frac{\sigma}{\varphi}\bm{g}(p))\notag\\
&\quad+\frac{\sigma}{8}\|\bar{\bm{g}}(p)\|^{2}\notag\\
&=V(p)-\epsilon_{8}\bm{x}^{T}(p)\bm{K}\bm{x}(p)+\epsilon_{2}\|\bm{e}(p)\|^{2}\notag\\
&\quad-\epsilon_{6}(\bm{u}(p)+\frac{\sigma}{\varphi}\bm{g}(p))^{T}\bm{P}(\bm{u}(p)+\frac{\sigma}{\varphi}\bm{g}(p))+\frac{\sigma}{8}\|\bar{\bm{g}}(p)\|^{2},
\end{align}
where the second inequality holds due to (\ref{LipIne}), (\ref{Lemma1Bb}), (\ref{Lemma1Bc}) and $\rho(\bm{K})=1$.

%
%
%

From the P--{\L} condition, similar to the (\ref{OptimalF2}), it can be calculated that
\begin{align}\label{OptimalF}
&n(f(\bar{x}(k+1))-f^{\star})\notag\\
&\leq n(f(\bar{x}(k))-f^{\star})-\frac{\sigma}{4}(1-2 \sigma L_{f})\|\bar{\bm{g}}(k)\|^{2}\notag\\
&\quad+\frac{\sigma}{2} L_{f}^{2}\bm{x}^{T}(k)\bm{K}\bm{x}(k)-\frac{\sigma}{2}\nu n(f(\bar{x}(k))-f^{\star}).\notag\\
\end{align}

From the Cauchy--Schwarz inequality, it is known that
\begin{align}\label{Vkl}
&V(k)\notag\\
&=\bm{x}^{T}(p)\bm{K}\bm{x}(p)+(\bm{u}(p)+\frac{\sigma}{\varphi}\bm{g}(p))^{T}(\frac{\varphi+\xi}{\varphi}\bm{P})(\bm{u}(p)\notag\\
&\quad+\frac{\sigma}{\varphi}\bm{g}(p))+2\bm{x}^{T}(p)\bm{K}\bm{P}(\bm{u}(p)+\frac{\sigma}{\varphi}\bm{g}(p))\notag\\
&\leq\bm{x}^{T}(p)\bm{K}\bm{x}(p)+(\bm{u}(p)+\frac{\sigma}{\varphi}\bm{g}(p))^{T}(\frac{\varphi+\xi}{\varphi}\bm{P})(\bm{u}(p)\notag\\
&\quad+\frac{\sigma}{\varphi}\bm{g}(p))+\frac{\varphi}{\xi\underline{\rho}(L)}\bm{x}^{T}(p)\bm{K}\bm{x}(p)+\frac{\xi}{\varphi}(\bm{u}(p)+\frac{\sigma}{\varphi}\bm{g}(p))^{T}\notag\\
&\quad\times\bm{P}(\bm{u}(p)+\frac{\sigma}{\varphi}\bm{g}(p))\notag\\
&\leq\frac{\xi\underline{\rho}(L)+\varphi}{\xi\underline{\rho}(L)}\bm{x}^{T}(p)\bm{K}\bm{x}(p)+(1+\frac{2\xi}{\varphi})(\bm{u}(p)+\frac{\sigma}{\varphi}\bm{g}(p))^{T}\notag\\
&\quad\times\bm{P}(\bm{u}(p)+\frac{\sigma}{\varphi}\bm{g}(p))\notag\\
&\leq\epsilon_{5}(\bm{x}^{T}(p)\bm{K}\bm{x}(p)+(\bm{u}(p)+\frac{\sigma}{\varphi}\bm{g}(p))^{T}\bm{P}(\bm{u}(p)+\frac{\sigma}{\varphi}\bm{g}(p))),
\end{align}
and
\begin{align}\label{Vkg}
&V(k)\notag\\
&=\bm{x}^{T}(p)\bm{K}\bm{x}(p)+(\bm{u}(p)+\frac{\sigma}{\varphi}\bm{g}(p))^{T}(\frac{\varphi+\xi}{\varphi}\bm{P})(\bm{u}(p)\notag\\
&\quad+\frac{\sigma}{\varphi}\bm{g}(p))+2\bm{x}^{T}(p)\bm{K}\bm{P}(\bm{u}(p)+\frac{\sigma}{\varphi}\bm{g}(p))\notag\\
&\geq\bm{x}^{T}(p)\bm{K}\bm{x}(p)+(\bm{u}(p)+\frac{\sigma}{\varphi}\bm{g}(p))^{T}(\frac{\varphi+\xi}{\varphi}\bm{P})(\bm{u}(p)\notag\\
&\quad+\frac{\sigma}{\varphi}\bm{g}(p))-\frac{\varphi}{\xi\underline{\rho}(L)}\bm{x}^{T}(p)\bm{K}\bm{x}(p)-\frac{\xi}{\varphi}(\bm{u}(p)+\frac{\sigma}{\varphi}\bm{g}(p))^{T}\notag\\
&\quad\times\bm{P}(\bm{u}(p)+\frac{\sigma}{\varphi}\bm{g}(p))\notag\\
&\geq\frac{\xi\underline{\rho}(L)-\varphi}{\xi\underline{\rho}(L)}\bm{x}^{T}(p)\bm{K}\bm{x}(p)+(\bm{u}(p)+\frac{\sigma}{\varphi}\bm{g}(p))^{T}\bm{P}(\bm{u}(p)\notag\\
&\quad+\frac{\sigma}{\varphi}\bm{g}(p))\notag\\
&\geq\epsilon_{10}(\bm{x}^{T}(p)\bm{K}\bm{x}(p)+(\bm{u}(p)+\frac{\sigma}{\varphi}\bm{g}(p))^{T}\bm{P}(\bm{u}(p)\notag\\
&\quad+\frac{\sigma}{\varphi}\bm{g}(p))).
\end{align}

From $\xi\geq\frac{5\varphi}{\underline{\rho}(L)}$, we know that $\frac{\xi\underline{\rho}(L)-\varphi}{\xi\underline{\rho}(L)}\geq\frac{4\varphi}{\xi\underline{\rho}(L)}>0$. Then, consider the following function:
\begin{equation}\label{Wk}
W(k)=V(k)+n(f(\bar{x}(k))-f^{\star}).
\end{equation}
Based on (\ref{Vkg}), we have
\begin{align}\label{Wkg}
&W(k)\notag\\
&\geq\epsilon_{10}(\bm{x}^{T}(p)\bm{K}\bm{x}(p)+(\bm{u}(p)+\frac{\sigma}{\varphi}\bm{g}(p))^{T}\bm{P}(\bm{u}(p)\notag\\
&\quad+\frac{\sigma}{\varphi}\bm{g}(p)))+n(f(\bar{x}_{k})-f^{\star})\notag\\
&\geq\epsilon_{10}(\bm{x}^{T}(p)\bm{K}\bm{x}(p)+(\bm{u}(p)+\frac{\sigma}{\varphi}\bm{g}(p))^{T}\bm{P}(\bm{u}(p)\notag\\
&\quad+\frac{\sigma}{\varphi}\bm{g}(p))+n(f(\bar{x}_{k})-f^{\star}))\notag\\
&=\epsilon_{10}\hat{W}(k)\geq0,
\end{align}
where $\hat{W}(k)=\bm{x}^{T}(p)\bm{K}\bm{x}(p)+(\bm{u}(p)+\frac{\sigma}{\varphi}\bm{g}(p))^{T}\bm{P}(\bm{u}(p)+\frac{\sigma}{\varphi}\bm{g}(p))+n(f(\bar{x}_{k})-f^{\star})$, and the second inequality holds due to $\epsilon_{10}\leq1$.
From (\ref{Vkl}), we have
\begin{align}\label{Wkl}
&W(k)\notag\\
&\leq\epsilon_{5}(\bm{x}^{T}(p)\bm{K}\bm{x}(p)+(\bm{u}(p)+\frac{\sigma}{\varphi}\bm{g}(p))^{T}\bm{P}(\bm{u}(p)\notag\\
&\quad+\frac{\sigma}{\varphi}\bm{g}(p)))+n(f(\bar{x}_{k})-f^{\star})\notag\\
&\leq\epsilon_{5}(\bm{x}^{T}(p)\bm{K}\bm{x}(p)+(\bm{u}(p)+\frac{\sigma}{\varphi}\bm{g}(p))^{T}\bm{P}(\bm{u}(p)\notag\\
&\quad+\frac{\sigma}{\varphi}\bm{g}(p))+n(f(\bar{x}_{k})-f^{\star}))\notag\\
&=\epsilon_{5}\hat{W}(k),
\end{align}
where the second inequality holds due to $\epsilon_{5}\geq1$. Then, \LX{for $\forall k=0,\dots,p$,} based on (\ref{Vk1l}) and (\ref{OptimalF}), it can be obtained that
\begin{align}\label{Wkk}
&W(k+1)\notag\\
&\leq W(k)-\epsilon_{8}\bm{x}^{T}(k)\bm{K}\bm{x}(k)+\epsilon_{2}\|\bm{e}(k)\|^{2}\notag\\
&\quad-\epsilon_{6}(\bm{u}(k)+\frac{\sigma}{\varphi}\bm{g}(k))^{T}\bm{P}(\bm{u}(k)+\frac{\sigma}{\varphi}\bm{g}(k))+\frac{\sigma}{8}\|\bar{\bm{g}}(k)\|^{2}\notag\\
&\quad-\frac{\sigma}{4}(1-2 \sigma L_{f})\|\bar{\bm{g}}(k)\|^{2}+\frac{\sigma}{2} L_{f}^{2}\bm{x}^{T}(k)\bm{K}\bm{x}(k)\notag\\
&\quad-\frac{\sigma}{2}\nu n(f(\bar{x}(k))-f^{\star})\notag\\
&= W(k)-\epsilon_{7}\bm{x}^{T}(k)\bm{K}\bm{x}(k)+\epsilon_{2}\|\bm{e}(k)\|^{2}\notag\\
&\quad-\epsilon_{6}(\bm{u}(k)+\frac{\sigma}{\varphi}\bm{g}(k))^{T}\bm{P}(\bm{u}(k)+\frac{\sigma}{\varphi}\bm{g}(k))\notag\\
&\quad-\frac{\sigma}{8}(1-4\sigma L_{f})\|\bar{\bm{g}}(k)\|^{2}-\frac{\sigma}{2}\nu n(f(\bar{x}(k))-f^{\star})\notag\\
&\leq W(k)-\epsilon_{4}(\bm{x}^{T}(k)\bm{K}\bm{x}(k)+(\bm{u}(k)+\frac{\sigma}{\varphi}\bm{g}(k))^{T}\bm{P}\notag\\
&\quad\times(\bm{u}(k)+\frac{\sigma}{\varphi}\bm{g}(k))+n(f(\bar{x}(k))-f^{\star}))+\epsilon_{2}\|\bm{e}(k)\|^{2}\notag\\
&\quad-\frac{\sigma}{8}(1-4\sigma L_{f})\|\bar{\bm{g}}(k)\|^{2}\notag\\
&= W(k)-\epsilon_{4}\hat{W}(k)+\epsilon_{2}\|\bm{e}(k)\|^{2}-\frac{\sigma}{8}(1-4\sigma L_{f})\|\bar{\bm{g}}(k)\|^{2}.
\end{align}


Denote
\begin{equation}\label{alpha1}
\alpha_{1}=\varphi-\frac{8\sigma}{\underline{\rho}(L)}-\frac{6\sigma^{2}\varphi^{2}L_{f}^{2}(\xi+\varphi)^{2}}{\varphi^{5}}-\frac{3\sigma L_{f}^{2}}{\underline{\rho}(L)},
\end{equation}
\begin{align}\label{alpha2}
\alpha_{2}&=\varphi^{2}\rho(L)+\frac{2\sigma^{2}}{\underline{\rho}(L)}+2\varphi^{2}\rho(L)\notag\\
&\quad+3\varphi^{2}L_{f}^{2}(\frac{\sigma^{2}(\xi+\varphi)}{\varphi^{3}}+\frac{3}{2}\rho(L)).
\end{align}

From (\ref{alpha1}), we have
\begin{align}\label{theta5}
\alpha_{1}&\geq\sigma(\kappa_{2}-\frac{6L_{f}^{2}(\kappa_{1}+1)^{2}}{\kappa_{2}}-\frac{3L_{f}^{2}+8}{\underline{\rho}(L)})\notag\\
&>\sigma(\kappa_{2}-1-\frac{3L_{f}^{2}+8}{\underline{\rho}(L)})\geq\sigma\varepsilon>\sigma^{2}\eta_{1},
\end{align}
where the first inequality holds due to $\xi\leq\kappa_{1}\varphi$ and $\varphi\geq\sigma\kappa_{2}$; the second inequality holds due to $\kappa_{2}\geq6L_{f}^{2}(\kappa_{1}+1)^{2}$; the third inequality holds due to $\kappa_{2}\geq1+\frac{3L_{f}^{2}+8}{\underline{\rho}(L)}$ and $0<\varepsilon\leq\kappa_{2}-1-\frac{3L_{f}^{2}+8}{\underline{\rho}(L)}$; the last inequality holds due to \blue{$\sigma<\frac{\varepsilon}{\eta_{1}}$}.

From (\ref{alpha2}), we have
\begin{align}\label{theta4}
\alpha_{2}&\leq\sigma^{2}(\kappa_{3}^{2}\rho(L)+\frac{2}{\underline{\rho}(L)}+2\kappa_{3}^{2}\rho(L)\notag\\
&\quad+3\kappa_{3}^{2}L_{f}^{2}(\frac{\kappa_{1}+1}{\kappa_{2}^{2}}+\frac{3}{2}\rho(L)))=\sigma^{2}\eta_{1},
\end{align}
where the inequality holds due to $\xi\leq\kappa_{1}\varphi$ and $\varphi\in[\sigma\kappa_{2},\sigma\kappa_{3}]$.

From Proposition \ref{prop}, we know that $\epsilon_{6}=\alpha_{1}-\alpha_{2}$. Then from (\ref{theta5}) and (\ref{theta4}), it can be obtained that
\blue{\begin{equation}\label{epsilon3}
\epsilon_{6}=\alpha_{1}-\alpha_{2}>0.
\end{equation}}

Denote
\begin{align}\label{alpha3}
\alpha_{3}&=\xi\underline{\rho}(L)-\frac{9\varphi}{2}-\sigma-\frac{6\sigma^{2}\xi^{2}L_{f}^{2}(\xi+\varphi)^{2}}{\varphi^{5}}\rho(L)\notag\\
&\quad-\frac{3\sigma L_{f}^{2}\xi^{2}}{\varphi^{2}},
\end{align}
\begin{align}\label{alpha4}
\alpha_{4}&=4\xi^{2}\rho^{2}(L)+2(\varphi(\xi+\varphi)\rho(L)+\sigma^{2}+\varphi^{2})\notag\\
&\quad+3\xi^{2}L_{f}^{2}(\frac{\sigma^{2}(\xi+\varphi)}{\varphi^{3}}\rho(L)+\frac{3}{2}\rho^{2}(L)).
\end{align}

From (\ref{alpha3}), we have
\begin{align}\label{theta1}
\alpha_{3}&\geq\sigma(\frac{\kappa_{2}}{2}-1-\frac{6L_{f}^{2}\kappa_{1}^{2}(\kappa_{1}+1)^{2}}{\kappa_{2}}\rho(L)-3L_{f}^{2}\kappa_{1}^{2})\notag\\
&>\sigma(\frac{\kappa_{2}}{2}-2-3L_{f}^{2}\kappa_{1}^{2})\geq \sigma(\varepsilon+\frac{1}{2}L_{f}^{2})>\sigma^{2}\eta_{2},
\end{align}
where the first inequality holds due to $\xi\in[\frac{5}{\underline{\rho}(L)}\varphi,\kappa_{1}\varphi]$ and $\varphi\geq\sigma\kappa_{2}$; the second inequality holds due to $\kappa_{2}\geq6L_{f}^{2}(\kappa_{1}+1)^{2}\kappa_{1}^{2}\rho(L)$; the third inequality holds due to $\kappa_{2}>4+6L_{f}^{2}\kappa_{1}^{2}+L_{f}^{2}$
and $0<\varepsilon\leq\frac{\kappa_{2}}{2}-2-3L_{f}^{2}\kappa_{1}^{2}-\frac{1}{2}L_{f}^{2}$; the last inequality holds due to \blue{$\sigma<\frac{\varepsilon}{\eta_{2}}$}.

From (\ref{alpha4}), we have
\begin{align}\label{theta2}
\alpha_{4}&\leq\sigma^{2}(4\kappa_{1}^{2}\kappa_{3}^{2}\rho^{2}(L)+2(\kappa_{3}^{2}(\kappa_{1}+1)\rho(L)+1+\kappa_{3}^{2})\notag\\
&\quad+3\kappa_{1}^{2}L_{f}^{2}((\kappa_{1}+1)\rho(L)+\frac{3}{2}\kappa_{3}^{2}\rho^{2}(L)))=\sigma^{2}\eta_{2},
\end{align}
where the inequality holds due to $\xi\leq\kappa_{1}\varphi$ and $\varphi\in[\sigma\kappa_{2},\sigma\kappa_{3}]$.

From Proposition \ref{prop}, we know that $\epsilon_{8}=\alpha_{3}-\alpha_{4}$. Then from (\ref{theta1}), (\ref{theta2}), it can be obtained that
\blue{\begin{equation}\label{epsilon5}
\epsilon_{8}=\alpha_{3}-\alpha_{4}>0.
\end{equation}}

Then from (\ref{epsilon5}), we have
\begin{align}\label{epsilon4}
\epsilon_{7}=\alpha_{3}-\alpha_{4}-\frac{\sigma}{2}L_{f}^{2}>\sigma(\varepsilon+\frac{1}{2}L_{f}^{2})-\alpha_{4}-\frac{\sigma}{2}L_{f}^{2}>0,
\end{align}
where the first inequality holds due to (\ref{theta1}); the second inequality holds due to (\ref{theta2}); \blue{the last inequality holds due to (\ref{theta2}) and $\sigma<\frac{\varepsilon}{\eta_{2}}$.}

Then, based on (\ref{epsilon3}), (\ref{epsilon4}), $\epsilon_{5}\geq1$ and $\sigma<\frac{2}{\nu}$, we have
\begin{equation}\label{epsilon6}
\epsilon_{3}=1-\frac{\epsilon_{4}}{\epsilon_{5}}=1-\frac{\min\{\epsilon_{6},\epsilon_{7},\frac{\sigma}{2}\nu\}}{\epsilon_{5}}\geq1-\frac{\sigma\nu}{2}>0.
\end{equation}
From (\ref{epsilon3}), (\ref{epsilon4}) and (\ref{epsilon6}), it can be guaranteed that $\epsilon_{3}\in(0,1)$.
Then, based on (\ref{Wkl}) and (\ref{Wkk}), we have
\begin{align}\label{VP1Cauuchy}
&W(p+1)\notag\\
&\leq\epsilon_{3}W(p)+\epsilon_{2}\|\bm{e}(p)\|^{2}-\frac{\sigma}{8}(1-4\sigma L_{f})\|\bar{\bm{g}}(p)\|^{2}\notag\\
&\leq\epsilon_{3}^{p+1}W(0)+\epsilon_{2}\sum_{\tau=0}^{p}\epsilon_{3}^{p-\tau}\|\bm{e}(\tau)\|^{2}\notag\\
&=\epsilon_{3}^{p+1}W(0)+\epsilon_{2}\sum_{\tau=0}^{p}s^{2}(\tau)\epsilon_{3}^{p-\tau}\|\frac{\bm{e}(\tau)}{s(\tau)}\|^{2}\notag\\
&\leq\epsilon_{3}^{p+1}W(0)+nm\epsilon_{2}s^{2}(p)\sum_{\tau=0}^{p}(\frac{\epsilon_{3}}{\mu^{2}})^{p-\tau}\frac{1}{4\mu^{2}}\notag\\
&=\epsilon_{3}^{p+1}W(0)+nm\epsilon_{2}\frac{s^{2}(p)}{4\mu^{2}}(\frac{\mu^{2}}{\mu^{2}-\epsilon_{3}})(1-(\frac{\epsilon_{3}}{\mu^{2}})^{p+1}),
\end{align}
where the second inequality holds due to $\sigma\leq\frac{1}{4L_{f}}$; the third inequality holds due to $\|x\|^{2}\leq nm\|x\|^{2}_{\infty}$ and (\ref{EpsilonK}).

From $\mu\in(\sqrt{\epsilon_{3}},1)$, it can be guaranteed that $\mu^{2}>\epsilon_{3}$. From $s(0)\geq\sqrt{\frac{4\mu^{2}(\mu^{2}-\epsilon_{3})W(0)}{\epsilon_{2}nm}}$ and $\mathcal{K}\geq\epsilon_{1}\sqrt{\frac{\epsilon_{2}nm}{4\mu^{2}(\mu^{2}-\epsilon_{3})}}+\frac{(1+2\xi d)}{2\mu}-\frac{1}{2}$, it can be guaranteed that
\begin{align}\label{thetaP1b}
&\|\theta(p+1)\|_{\infty}\notag\\
&\leq\frac{1}{s(p+1)}\epsilon_{1}(V(p+1))^{\frac{1}{2}}+\frac{(1+2\xi d)}{2\mu}\notag\\
&\leq\frac{1}{s(p+1)}\epsilon_{1}(W(p+1))^{\frac{1}{2}}+\frac{(1+2\xi d)}{2\mu}\notag\\
&\leq\frac{1}{s(p+1)}\epsilon_{1}(nm\epsilon_{2}\frac{s^{2}(p)}{4\mu^{2}}(\frac{\mu^{2}}{\mu^{2}-\epsilon_{3}})(1-(\frac{\epsilon_{3}}{\mu^{2}})^{p+1})\notag\\
&\quad+\epsilon_{3}^{p+1}W(0))^{\frac{1}{2}}+\frac{(1+2\xi d)}{2\mu}\notag\\
&=\frac{(1+2\xi d)}{2\mu}+\epsilon_{1}(\frac{W(0)}{s^{2}(0)}(\frac{\epsilon_{3}}{\mu^{2}})^{p+1}\notag\\
&\quad+\frac{nm\epsilon_{2}}{4\mu^{2}(\mu^{2}-\epsilon_{3})}(1-(\frac{\epsilon_{3}}{\mu^{2}})^{p+1}))^{\frac{1}{2}}\notag\\
&\leq \LX{\Omega+\frac{1}{2}}\leq\mathcal{K}+\frac{1}{2},
\end{align}
where the first inequality holds due to (\ref{thetaP1a});
the second inequality holds due to $f(\bar{x}(k))-f^{\star}\geq0$; the third inequality holds due to (\ref{VP1Cauuchy}).

As a result, when $k=p+1$, the quantizer is also unsaturated. Therefore, by induction, we conclude that the quantizer is never saturated. This proof is complete.

\section{Proof of Theorem \ref{thm-1}}\label{Appendix-B}

From (\ref{VP1Cauuchy}) and $s(0)\geq\sqrt{\frac{4\mu^{2}(\mu^{2}-\epsilon_{3})W(0)}{\epsilon_{2}nm}}$, we have
\begin{align}\label{Wk1l}
&W(k+1)\notag\\
&\leq\epsilon_{3}^{k+1}W(0)+nm\epsilon_{2}\frac{s^{2}(k)}{4\mu^{2}}(\frac{\mu^{2}}{\mu^{2}-\epsilon_{3}})(1-(\frac{\epsilon_{3}}{\mu^{2}})^{k+1})\notag\\
&=\epsilon_{3}^{k+1}(W(0)-\frac{nm\epsilon_{2}s^{2}(0)}{4\mu^{2}(\mu^{2}-\epsilon_{3})})+\frac{nm\epsilon_{2}s^{2}(0)}{4\mu^{2}(\mu^{2}-\epsilon_{3})}\mu^{2(k+1)}\notag\\
&\leq\frac{nm\epsilon_{2}s^{2}(0)}{4\mu^{2}(\mu^{2}-\epsilon_{3})}\mu^{2(k+1)},
\end{align}

Then, from (\ref{Wkg}), we have
\begin{equation}\label{hatWk}
\bm{x}^{T}(p)\bm{K}\bm{x}(p)+n(f(\bar{x}_{k})-f^{\star})\leq\hat{W}(k)\leq\frac{W(k)}{\epsilon_{10}}.
\end{equation}
Hence, (\ref{Wk1l}) and (\ref{hatWk}) yield (\ref{Thm1A}). This proof is complete.

\section{Proof of Theorem \ref{thm-2}}\label{Appendix-C}


From $\xi\in[\frac{5}{\underline{\rho}(L)}\varphi,\kappa_{1}\varphi]$ and $\varphi\in[\sigma\kappa_{2},\sigma\kappa_{3}]$, it can be found that
\begin{align}\label{eta0}
\lim_{\sigma\rightarrow0}\epsilon_{1}\sqrt{\frac{\epsilon_{2}nm}{4\mu^{2}(\mu^{2}-\epsilon_{3})}}+\frac{(1+2\xi d)}{2\mu}\leq\frac{1}{2\mu},
\end{align}
where the inequality holds due to $\lim_{\sigma\rightarrow0}\epsilon_{2}=0$ and $\lim_{\sigma\rightarrow0}\xi=0$.

Then, for any given $\mathcal{K}\geq1$, there exists $\sigma^{\star}\in(0,\min\{\blue{\frac{\varepsilon}{\eta_{1}},\frac{\varepsilon}{\eta_{2}}},\\
\frac{2}{\nu},\frac{1}{4L_{f}}\})$ such that
\begin{align}\label{mu1}
&\lim_{\mu\rightarrow1}\epsilon_{1}\sqrt{\frac{\epsilon_{2}nm}{4\mu^{2}(\mu^{2}-\epsilon_{3})}}+\frac{(1+2\sigma^{\star}\xi d)}{2\mu}\leq\mathcal{K}+\frac{1}{2}.
\end{align}

Hence, there exists $\mu^{\star}\in(\sqrt{\epsilon_{3}},1)$ such that $\LX{\Omega\leq\mathcal{K}}$. Thus, $(\mu^{\star},\sigma^{\star})\in\bar{\Pi}$, and hence $\bar{\Pi}$ is nonempty, where $\bar{\Pi}$ is defined in Theorem \ref{thm-2}. The proof of the convergence result is similar to that of Theorem \ref{thm-1}.
\end{document}